\input amstex
\documentstyle{amsppt}
\input colordvi
\let\a=\alpha \let\b=\beta  \let\d=\delta \let\e=\varepsilon
 \let\g=\gamma \let\h=\eta \let\k=\kappa 
   
\let\r=\rho \let\s=\sigma  \let \th=\theta 
 \let\x=\xi 
\let\D=\Delta

    \let\X=\Xi

\def\data{\number\day/\ifcase\month\or gennaio \or febbraio \or marzo \or
aprile \or maggio \or giugno \or luglio \or agosto \or settembre
\or ottobre \or novembre \or dicembre \fi/\number\year;\,\the\time}
\def\dat{\number\day\,\ifcase\month\or {\rm gennaio} \or {\rm febbraio} \or
{\rm marzo}
\or {\rm aprile} \or {\rm maggio} \or {\rm giugno} \or {\rm luglio} \or {\rm
agosto}
\or {\rm settembre} \or {\rm ottobre} \or {\rm novembre} \or {\rm dicembre}
\fi\,\number\year}
\newcount\pgn \pgn=1
\def\foglio{\number\numsec:\number\pgn
\global\advance\pgn by 1}
\def\foglioa{A\number\numsec:\number\pgn
\global\advance\pgn by 1}

\global\newcount\numsec\global\newcount\numfor
\gdef\profonditastruttura{\dp\strutbox}
\def\senondefinito#1{\expandafter\ifx\csname#1\endcsname\relax}
\def\SIA #1,#2,#3 {\senondefinito{#1#2}%
\expandafter\xdef\csname #1#2\endcsname{#3}\else
\write16{\?? ma #1,#2 e' gia' stato definito !!!!} \fi}
\def\etichetta(#1){(\veroparagrafo.\veraformula)%
\SIA e,#1,(\veroparagrafo.\veraformula) %
\global\advance\numfor by 1%
\write15{\string\FU (#1){\equ(#1)}}%
}
\def\FU(#1)#2{\SIA fu,#1,#2 }
\def\etichettaa(#1){(A.\veraformula)%
\SIA e,#1,(A.\veraformula) %
\global\advance\numfor by 1%
\write15{\string\FU (#1){\equ(#1)}}%
}
\def\getichetta(#1){Fig. \verafigura
\SIA e,#1,{\verafigura} %
\global\advance\numfig by 1%
\write15{\string\FU (#1){\equ(#1)}}%
\write16{ Fig. \equ(#1) ha simbolo #1 }}
\newdimen\gwidth
\def\BOZZA{
\def\alato(##1){%
 {\vtop to \profonditastruttura{\baselineskip
 \profonditastruttura\vss
 \rlap{\kern-\hsize\kern-1.2truecm{$\scriptstyle##1$}}}}}
\def\galato(##1){\gwidth=\hsize \divide\gwidth by 2%
 {\vtop to \profonditastruttura{\baselineskip
 \profonditastruttura\vss
 \rlap{\kern-\gwidth\kern-1.2truecm{$\scriptstyle##1$}}}}}
}
\def\alato(#1){}
\def\galato(#1){}
\def\veroparagrafo{\number\numsec}\def\veraformula{\number\numfor}
\def\verafigura{\number\numfig}
\def\geq(#1){\getichetta(#1)\galato(#1)}
\def\Eq(#1){\eqno{\etichetta(#1)\alato(#1)}}
\def\eq(#1){\etichetta(#1)\alato(#1)}
\def\Eqa(#1){\eqno{\etichettaa(#1)\alato(#1)}}
\def\eqa(#1){\etichettaa(#1)\alato(#1)}
\def\eqv(#1){\senondefinito{fu#1} #1
\write16{#1 non e' (ancora) definito}%
\else\csname fu#1\endcsname\fi}
\def\equ(#1){\senondefinito{e#1}\eqv(#1)\else\csname e#1\endcsname\fi}
\def\include#1{
\openin13=#1.aux \ifeof13 \relax \else
\input #1.aux \closein13 \fi}
\openin14=\jobname.aux \ifeof14 \relax \else
\input \jobname.aux \closein14 \fi
\openout15=\jobname.aux

 
 \def\\{\noindent}

\let \pt=\partial

\def\tende#1{\vtop{\ialign{##\crcr\rightarrowfill\crcr
 \noalign{\kern-1pt\nointerlineskip}
 \hskip3.pt${\scriptstyle #1}$\hskip3.pt\crcr}}}
\def\otto{{\kern-1.truept\leftarrow\kern-5.truept\to\kern-1.truept}}

\def\mbox{\hbox}
\def\Tr{{\rm Tr}\,}

\def\={{\equiv}}

\def\n{\nabla}

\def\alf{{1\over2}}
\def\reset{}

\def\11{\hbox{l}\!\!\!\hbox{1}\,}

\def\sqr#1#2{{\vcenter{\vbox{\hrule height.#2pt
\hbox{\vrule width.#2pt height #1pt \kern#1pt
\vrule width.#2pt}
\hrule height.#2pt}}}}
\def\square{\mathchoice\sqr56\sqr56\sqr{2.1}3\sqr{1.5}3}
\def\qed{$\square$}

\def\picture #1 by #2 (#3){
 \vbox to #2{
 \hrule width #1 height 0pt depth 0pt
 \vfill
 \special{picture #3}
 }
 }
\def\scaledpicture #1 by #2 (#3 scaled #4){{
 \dimen0=#1 \dimen1=#2
 \divide\dimen0 by 1000 \multiply\dimen0 by #4
 \divide\dimen1 by 1000 \multiply\dimen1 by #4
 \picture \dimen0 by \dimen1 (#3 scaled #4)}
 }

\def\und(#1){$\underline{\hbox{#1}}$}

\def\xp(#1){\hbox{\rm e}^{#1}}

\def\noi{\noindent}
\def\({\Big(}
\def\){\Big)}
\def\ex{\hbox{\rm e}}


\def\Tr{\hbox{\rm Tr\,}}

\catcode`\X=12\catcode`\@=11
\def\n@wcount{\alloc@0\count\countdef\insc@unt}
\def\n@wwrite{\alloc@7\write\chardef\sixt@@n}
\def\n@wread{\alloc@6\read\chardef\sixt@@n}
\def\crossrefs#1{\ifx\alltgs#1\let\tr@ce=\alltgs\else\def\tr@ce{#1,}\fi
   \n@wwrite\cit@tionsout\openout\cit@tionsout=\jobname.cit 
   \write\cit@tionsout{\tr@ce}\expandafter\setfl@gs\tr@ce,}
\def\setfl@gs#1,{\def\@{#1}\ifx\@\empty\let\next=\relax
   \else\let\next=\setfl@gs\expandafter\xdef
   \csname#1tr@cetrue\endcsname{}\fi\next}
\newcount\sectno\sectno=0\newcount\subsectno\subsectno=0\def\r@s@t{\relax}
\def\resetall{\global\advance\sectno by 1\subsectno=0
   \gdef\firstpart{\number\sectno}\r@s@t}
\def\resetsub{\global\advance\subsectno by 1
   \gdef\firstpart{\number\sectno.\number\subsectno}\r@s@t}
\def\v@idline{\par}\def\firstpart{\number\sectno}
\def\l@c@l#1X{\firstpart.#1}\def\gl@b@l#1X{#1}\def\t@d@l#1X{{}}
\def\m@ketag#1#2{\expandafter\n@wcount\csname#2tagno\endcsname
     \csname#2tagno\endcsname=0\let\tail=\alltgs\xdef\alltgs{\tail#2,}%
   \ifx#1\l@c@l\let\tail=\r@s@t\xdef\r@s@t{\csname#2tagno\endcsname=0\tail}\fi
   \expandafter\gdef\csname#2cite\endcsname##1{\expandafter
     \ifx\csname#2tag##1\endcsname\relax?\else{\rm\csname#2tag##1\endcsname}\fi
     \expandafter\ifx\csname#2tr@cetrue\endcsname\relax\else
     \write\cit@tionsout{#2tag ##1 cited on page \folio.}\fi}%
   \expandafter\gdef\csname#2page\endcsname##1{\expandafter
     \ifx\csname#2page##1\endcsname\relax?\else\csname#2page##1\endcsname\fi
     \expandafter\ifx\csname#2tr@cetrue\endcsname\relax\else
     \write\cit@tionsout{#2tag ##1 cited on page \folio.}\fi}%
   \expandafter\gdef\csname#2tag\endcsname##1{\global\advance
     \csname#2tagno\endcsname by 1%
   \expandafter\ifx\csname#2check##1\endcsname\relax\else%
\fi
   \expandafter\xdef\csname#2check##1\endcsname{}%
   \expandafter\xdef\csname#2tag##1\endcsname
     {#1\number\csname#2tagno\endcsnameX}%
   \write\t@gsout{#2tag ##1 assigned number \csname#2tag##1\endcsname\space
      on page \number\count0.}%
   \csname#2tag##1\endcsname}}%
\def\m@kecs #1tag #2 assigned number #3 on page #4.%
    {\expandafter\gdef\csname#1tag#2\endcsname{#3}
    \expandafter\gdef\csname#1page#2\endcsname{#4}}
\def\re@der{\ifeof\t@gsin\let\next=\relax\else
    \read\t@gsin to\t@gline\ifx\t@gline\v@idline\else
    \expandafter\m@kecs \t@gline\fi\let \next=\re@der\fi\next}
\def\t@gs#1{\def\alltgs{}\m@ketag#1e\m@ketag#1s\m@ketag\t@d@l p
    \m@ketag\gl@b@l r \n@wread\t@gsin\openin\t@gsin=\jobname.tgs \re@der
    \closein\t@gsin\n@wwrite\t@gsout\openout\t@gsout=\jobname.tgs }
\outer\def\localtags{\t@gs\l@c@l}
\outer\def\globaltags{\t@gs\gl@b@l}
\outer\def\newlocaltag#1{\m@ketag\l@c@l{#1}}
\outer\def\newglobaltag#1{\m@ketag\gl@b@l{#1}}

\def\t@gsoff#1,{\def\@{#1}\ifx\@\empty\let\next=\relax\else\let\next=\t@gsoff
   \expandafter\gdef\csname#1cite\endcsname{\relax}
   \expandafter\gdef\csname#1page\endcsname##1{?}
   \expandafter\gdef\csname#1tag\endcsname{\relax}\fi\next}
\def\verbatimtags{\let\ift@gs=\iffalse\ifx\alltgs\relax\else
   \expandafter\t@gsoff\alltgs,\fi}
\def\(#1){\edef\dot@g{\ift@gs\if@lign(\noexpand\etag{#1})
    \else\ifmmode\noexpand\tag\noexpand\etag{#1}%
    \else{\rm(\noexpand\ecite{#1})}\fi\fi\else\if@lign{\rm(#1)}
    \else\ifmmode\noexpand\tag#1\else{\rm(#1)}\fi\fi\fi}\dot@g}
\let\ift@gs=\iftrue\let\if@lign=\iffalse
\let\n@weqalignno=\eqalignno
\def\eqalignno#1{\let\if@lign=\iftrue\n@weqalignno{#1}}
\catcode`\X=11 \catcode`\@=\active
\localtags
\nologo \magnification=1200 \hsize=16.5truecm
\vsize=24.truecm \topmatter


\title
On The Weak-Coupling Limit for Bosons and Fermions
\endtitle
\author
D. Benedetto$^*$
\footnote""{$^*$ \eightrm Dipartimento di Matematica
Universit\`a di Roma `La Sapienza', P.le A. Moro, 5 - 00185 Roma - Italy},
F. Castella$^+$
\footnote""{$^+$ \eightrm IRMAR, Universit\'e de Rennes 1, Campus de Beaulieu
- 35042 Rennes - Cedex - France
\quad\quad\quad\quad},
R. Esposito$^\#$
\footnote""{$^\#$ \eightrm Dipar\-timento di Matematica
pura ed applicata, Universit\`a di L'Aquila, Coppito -  67100 L'Aquila - Italy}
and
M. Pulvirenti$^*$
\endauthor
\bigskip
\bigskip
\abstract
In this paper we consider a large system of Bosons or Fermions. We
start with an initial datum which is compatible with the Bose-Einstein,
respectively Fermi-Dirac, statistics. We let the system of interacting particles
evolve in a weak-coupling regime. We show that, in the limit, 
and up to the second order in the potential,
the
perturbative expansion  expressing the value of the one-particle Wigner
function at time $t$, agrees with the analogous expansion for the solution to the
Uehling-Uhlenbeck equation.

This paper follows in spirit the companion work [\rcite{BCEP}], where the 
authors
investigated the weak-coupling limit for particles obeying the 
Maxwell-Boltzmann  statistics: here,
they proved a (much stronger) convergence result towards the solution
of the Boltzmann equation.
\endabstract
\date
3/10/2004
\enddate
\endtopmatter


\heading 1. Introduction \endheading
\numsec= 1
\reset
\numfor= 1


In 1933 Uehling and Uhlenbeck in Ref. [\rcite{UU}] proposed the following
kinetic equation, called U-U in the sequel, for the time 
evolution of the one-particle Wigner function $f(x,v;t)$
associated with a
large system of weakly interacting Bosons or Fermions
(see Ref. [\rcite{W}] for the definition of the Wigner function).
The U-U equation is
$$
\eqalign{
\pt_{t}f(x,v;t)+&v\cdot\n_{x}f(x,v;t) =
\int dv_{*}\int dv_{*}'\int dv'\;
W(v,v_{*}|v',v_{*}')
\cr& 
\Big\{f'f'_{*}(1+8\pi^3\th f)(1+8\pi^3\th f_{*})
-ff_{*}(1+8\pi^3\th f')(1+8\pi^3\th f_{*}')\Big\},}
\Eq(UU)
$$
where we use the standard short-hand notation
$$
f=f(x,v;t),
\quad f_{*}=f(x,v_{*};t),
\quad f'=f(x,v';t),
\quad f'_{*}=f(x,v'_{*};t).
$$
Here,  $W$ denotes the
transition kernel
$$
\eqalign{
W(v,v_{*}| v',v'_{*})  =
{1\over 8\pi^{2}}
&
\big[\widehat \phi(v'-v)+\th \, \widehat \phi(v'-v_{*}) \big]^{2}
\cr &
\d(v+v_{*}-v'-v'_{*})
\; \d\Big(\alf(v^{2}+v_{*}^{2}-v'\/^{2}-v'_{*}\/^{2})\Big).}
\Eq(W)
$$
Finally,
$$
\widehat \phi(k)= \int dx \; \ex^{-ik\cdot x}\phi(x)
\hskip4cm
\Eq(Fou)
$$
is the Fourier transform of the two-body interaction potential $\phi$, and
$\th=\pm 1$ for Bosons and Fermions respectively.

Note that the factors $8\pi^3$ do not appear in the original U-U equation
in Ref. [\rcite{UU}], because there, the distribution function is
normalized  in such a way that its integral on the velocity variable
equals the space density times $8\pi^3$.
At variance, in \equ(UU), $f$ is just the
standard Wigner function, whose integral on the velocities equals the space
density. Let us mention that equation \equ(UU) actually is
cubic (and not quartic)
in the unknown $f$: apparent quartic terms have vanishing contribution, as
shown by direct inspection.

Eq. \equ(UU) constitutes a natural modification of the usual
quantum Boltzmann
equation, in order to take into account statistics. In particular, there is a
$H$-functional
$$
\Cal H(f)=
\int dx\int dv \, \big\{f\log f -\th(1+8\pi^3\th f)\log (1+8\pi^3\th f)
\big\}
\Eq(H)
$$
driving the
system to the Bose-Einstein and Fermi-Dirac equilibrium distribution:
$$
M(v)={1\over \ex^{\b\mu+\b v^2/2} -8\pi^3\th}
\Eq(M)
$$
outside the Bose condensation region. Here $\b$ and $\mu$ denote the
inverse temperature and the chemical potential respectively.

Eq. \equ(UU) is largely studied (see for istance
[\rcite{B}] and [\rcite{J}] for physical consideration, 
and [\rcite{LW}], [\rcite{D}], [\rcite{LX1}], [\rcite{LX2}] ... for 
a more mathematically oriented analysis concerning the 
existence of solutions and asymptotic behavior),
so that it is 
certainly of great relevance to derive this equation from
the first principles, namely from the Schr\"odinger equation.

As clearly explained by H. Spohn in [\rcite{S}], Eq \equ(UU) is indeed
expected to hold in the so called weak-coupling limit, which consists in
scaling space, time and the potential according to
$$
x\to \e x, \quad
t\to \e t, \quad
\phi\to \sqrt \e \phi,
\Eq(sc)
$$
where $\e$ is  a small positive parameter.

A slightly different limit, usually called van Hove limit, scales $t$ and
$\phi$ as in \equ(sc) but leaves the microscopic space scale unchanged.
Eq. \equ(UU) cannot be derived in the van Hove limit in general but, in
case of translationally invariant states, we expect to achieve the
homogeneous version of the U-U equation { (for a large system)}. In fact
Hugenholtz [\rcite{H}] proved formally that this happens. Later on Ho and
Landau [\rcite{HL}] proved that the homogeneous U-U equation holds
rigorously up to the second order expansion in the potential. These
approaches, as well as the recent contribution by Erd\"os et al.
[\rcite{ESY}] (where the quantum analog of the Boltzmann's {\it
Stosszahlansatz} is formulated), are based on the commutator expansion of
the time evolution of the observables of the CCR and CAR algebras.

In the present paper we approach the problem from a different viewpoint.
We start from the time evolution of a $N$ particle quantum system in terms
of the Wigner formalism. Here the statistics enters only through the choice of
the admissible states we take as initial conditions. Such states,
called quasi-free, must describe free Bosons and Fermions, so that they
cannot have any other correlations but those arising from the statistics.
Therefore the first step is to characterize quasi-free states (see for
example Ref. [\rcite{BR}]) in terms of the Wigner functions. Then we apply
the dynamics (in terms of the usual hierarchy) and represent the solution
as a perturbative expansion. The truncation of this expansion up to the
second order in the potential is shown to converge to the expansion
associated to the U-U equation, up to the first order in the scattering
cross section.

In other words we recover the result in Ref. [\rcite{HL}] with the
following main differences.  First we exploit the weak-coupling
limit, so that we can deal with states which are not necessarily
translationally invariant. Second, we work directly in terms of
the Wigner formalism, in the same spirit of the Balescu book (see Ref.
[\rcite{B}]). In doing so, we also follow a previous work 
[\rcite{BCEP}] by the
authors for the Maxwell-Boltzmann statistics.
Hence the present work shows how the statistics can be handled in this
formalism. Note in passing that the case of the Maxwell-Boltzmann 
statistics allows for a much stronger (but still partial)
convergence result than the one 
presented here, see [\rcite{BCEP}]. Note finally that the present formalism
also allows to handle the low-density limit, see [\rcite{BCEP2}],
see also the last section of this text.

It is also important to mention that a full rigorous derivation of the U-U
equation (but also of the usual Boltzmann equation arising for the
Maxwell-Boltzmann statistics) is still far beyond the present techniques
and those of the previous references.

The plan of the paper is the following. In the next section we describe
the particle system. In Section 3 we establish the result. The rest of the
paper is devoted to the proofs.

\bigskip\bigskip

\heading 2. The particle system.\endheading
\numsec= 2
\reset
\numfor= 1
\bigskip

We consider a Quantum particle system in $\Bbb R^{3}$. Let
$$
\frak H =
\bigoplus_{n\ge 0}L_{2}(\Bbb R^{3})^{n}
:= \bigoplus_{n\ge 0}\frak H_{n},
\Eq(2.1)
$$
be the Fock space. A state of 
the system is a self-adjoint,
positive trace class operator acting on $\frak H$:
$$
\s=\bigoplus_{n\ge 0}\s_{n}.
\Eq(2.2)
$$
We assume
$$
\hbox{\rm Tr\,} \s=1.
\Eq(2.3)
$$
The operator $N$, number of particles, is the multiplication by $n$ on
$\frak H_{n}$ and hence
$$
\langle N\rangle=
\sum_{n \ge 0}n \,\hbox{\rm Tr\,} \s_{n},
\Eq(2.4)
$$
where the left hand side is the average number of particles in the state $\s$.
If $\s_{n}(X_{n};Y_{n})$ is the kernel of $\s_{n}$,  the Reduced
Density Matrices (RDM) are defined by:
$$
\r_{n}(X_{n};Y_{n})=
\sum_{m\ge 0}\frac{(n+m)!}{m!}
\int \s_{n+m}(X_{n},Z_{m};Y_{n},Z_{m}) \, dZ_{m}.
\Eq(2.5)
$$ 
Here $X_n=(x_1,\dots,x_n)$, $x_i\in \Bbb R^3$ denotes the $n$-particle
configuration.
Note that
$$
\eqalign{
\Tr \r_{n} & =
\int dZ_n \; \r_{n}(Z_{n};Z_{n})
\cr&
=\sum_{m \ge n}m(m-1)\dots(m-n+1) 
\Tr \s_{m}=\langle N(N-1)\dots (N-n+1)\rangle,
}\Eq(2.6)
$$
and hence the RDM 
are equivalent
to the classical correlation functions.

The Hamiltonian of the system is the self-adjoint
operator acting
on $\frak H$ given by
$$
H=\bigoplus_{n=1}^{\infty}H_{n},
\Eq(2.7)
$$
where
$$
H_{n}=
-{1\over 2}\sum_{j=1}^{n}\D_{x_{j}}
+\sum_{1\le i<j\le n}\phi (x_{i}-x_{j}),
\Eq(2.8)
$$
and the potential $\phi$ is a smooth two-body interaction.
Here, $\hbar $ as 
well as the mass
of the particles are normalized to unity.

Under these circumstances, the time evolved state is given by
the usual
$$
\s(t)=\ex^{-iHt}\s\ex^{iHt}.
\Eq(2.9)
$$

Now, quantum statistics is taken into account by suitable properties of
the physically relevant states. Namely, for the Maxwell-Boltzmann (M-B)
statistics we require symmetry of
$\r_{n}(x_{1},\dots,x_{n};y_{1}\dots,y_{n})$ in the exchange of particle
names. For the Bose-Einstein (B-E) and Fermi-Dirac (F-D) statistics we
require additionally
$$
\r_{n}(x_{1},\dots,x_{n};y_{1}\dots,y_{n})
=
\th^{s(\pi)}\r_{n}(x_{1},\dots,x_{n};y_{\pi(1)}\dots,y_{\pi(n)}),
\Eq(stat)
$$
where $\pi\in \Cal P_{n}$ is a
permutation of $n$ elements, and $s(\pi)=0$ if the permutation is even,
$s(\pi)=1$ if it is odd.

Alternatively, the quantum statistics is automatically taken into account by
considering states on the algebra generated by the annihilation
and creation operators $a(x)$ and $a^{\dagger}(x)$ (with the
commutation and anti-commutation relations according to the B-E
and F-D statistics respectively). Then the RDM are defined as
$$
\Tr\big[ \s a^{\dagger}(x_{n})\dots a^{\dagger}(x_{1})
a(y_{1})\dots a(y_{n})\big]= 
\r_{n}(x_{1},\dots,x_{n};y_{1}\dots,y_{n}).
\Eq(2.11)
$$
However we do not use 
here the second quantization formalism.

Given a state $\s$, we define the Wigner transform [\rcite{W}] by
$$
W_{n}(X_{n};V_{n}):=
{1\over (2\pi)^{3n}}\int d Y_{n} \;
\ex^{iY_{n}\cdot V_{n}}
\s_{n}\Big(X_{n}-{1\over 2} Y_{n};X_{n}+{1\over 2} Y_{n}\Big).
\Eq(wig)
$$ 
Therefore the analogous of the
classical correlation functions are the $j$-particle Wigner functions
defined through
$$
F_{j}(X_{j};V_{j})=
\sum_{n\ge 0}\frac{(n+j)!}{n!}
\int dX_{n}\int dV_{n} \;
W_{j+n}(X_{j},X_{n};V_{j},V_{n}).
\Eq(wigj)
$$  
Note that the
$F_{j}$'s are the Wigner transforms of the RDM $\r_{j}$, as one can easily 
check.

Due to the dynamics imposed by \equ(2.9), it is a standard computation to 
check that the Wigner 
function $W_{n}$ evolves according to the Wigner-Liouville equation 
$$
\pt_{t} W_{n}+\sum_{i=1}^{n}v_{i}\cdot \n_{x_{i}}W_{n}
=T_{n}W_{n}.
\Eq(wigheq)
$$ 
As a consequence,
the $j$-particle Wigner functions $F_{j}$'s satisfy the 
associated hierarchy
$$
\pt_{t} F_{j}+\sum_{i=1}^{j}v_{i}\cdot \n_{x_{i}}F_{j}
=T_{j}F_{j}+C_{j+1}F_{j+1},
\Eq(wigh)
$$
where $T_j$ and $C_{j+1}$ will be defined below after
Eq.\equ(2.21), and the index $j$ takes any value between $1$ and $N$.
Equations \equ(wigh) are analogous to the usual BBGKY hierarchy
for the classical systems and are derived in a similar manner. Note that
by the definition of the RDM the coefficient in front of $C_{j+1}$ is one
instead of $N-j$.

\bigskip

We now want to analyze \equ(wigh) in the
weak-coupling regime \equ(sc).
Therefore, we set
$$
f^{\e}_{j}(X_{j};V_{j};t):=
F_{j}(\e^{-1}X_{j};V_{j};\e^{-1}t),
\Eq(scal)
$$
where $\e>0$ is a small parameter, and we scale the potential 
as well, by setting
$$
\phi\to \sqrt\e \phi.
\Eq(2.16)
$$
The resulting, scaled, equation is
$$
\pt_{t}f^{\e}_{j}+\sum_{i=1}^{j}v_{i}\cdot \n_{x_{i}}f^{\e}_{j}
=
{1\over\sqrt\e}T^{\e}_{j} f^{\e}_{j}
+
{1\over\sqrt\e}\e^{-3}C^{\e}_{j+1}f^{\e}_{j+1},
\Eq(rp)
$$
where 
$$
(T^\e_j f_j^\e)(X_j;V_j)
= \sum_{0<k<\ell\le j} (T^\e_{k,l} f_j^\e)(X_j;V_j),
\Eq(2.18)
$$ 
and the $T^\e_{k,l}$'s are defined as follows: if $j=1$, 
we simply have $T_1^\e=0$; otherwise, 
$$
\eqalign{
(T^\e_{k,l} f_j^\e)(X_j;V_j)=&
-i\sum_{\s=\pm 1}\s\int  {dh \over (2\pi)^3}\,
\widehat\phi(h) \, \ex^{i{h\over \e}\cdot(x_k-x_\ell)}\cr&
f_j^\e\Big(x_1,\dots,x_j;
v_1,\dots, v_k-\s{h\over 2},\dots,v_\ell+\s{h\over 2},\dots, v_j\Big).}
\Eq(2.19)
$$
On the other hand, the $C^\e_{j+1}$ in \equ(rp) is computed as:
$$
(C^\e_{j+1} f_{j+1}^\e)(X_j;V_j)=
\sum_{k=1}^j (C^\e_{k,j+1} f_{j+1}^\e)(X_j;V_j),
\Eq(2.20)
$$
where
$$
\eqalign{
(C^\e_{k,j+1} f_{j+1}^\e)(X_j;V_j)=&
-i\sum_{\s=\pm 1}\s\int {dh \over (2\pi)^3}
\int dx_{j+1}\int dv_{j+1} \; 
\widehat\phi(h) \, \ex^{i{h\over \e}\cdot(x_k-x_{j+1})}\cr&
f_{j+1}^\e\Big(x_1,\dots,x_{j+1};
v_1,\dots, v_k-\s{h\over 2},\dots,v_{j+1}+\s{h\over 2}\Big).}
\Eq(2.21)
$$
Note that $T_j$ and $C_{j+1}$ are $T_j^\e$ and $C_{j+1}^\e$ for $\e=1$.
Last, we fix an initial condition sequence
$$
\{f_{0,j}^{\e}\}_{j=1}^{\infty}
\Eq(is)
$$
according to the quantum statistics, and perform the limit $\e\to 0$
in the resulting system.

\bigskip

\noi
$\underline{\hbox{\rm Remark}}$: Since
$$
\int f_{0,1}^{\e}(x,v)dxdv=\e^{3}\langle N\rangle,
\Eq(2.24)
$$
requiring $\|f_{0,1}^{\e}\|_{L_{1}}=O(1)$, implies
$\langle N\rangle=O(\e^{-3})$. In other 
words we are working in the Grand-canonical formalism and the density is
automatically rescaled.
\hfill\qed

\bigskip

In the following we shall fix $ f_{0,1}^{\e}$ to be a given (independent
of $\e$) probability density $ f_{0}$. This means that its inverse Wigner
transform
$$
\r^{\e}(x,y)=
\int dv \,
\ex^{i {x-y\over \e} \cdot v} \,
f_{0}\Big({x+y\over 2}, v\Big),
\Eq(r)
$$
namely the one-particle rescaled RDM, is a superposition of WKB states.

We now make assumptions on the initial state to take into account the
statistics.
For the M-B statistics a suitable initial sequence can be
chosen completely factorized, e.g.
$$
f^{\e}_{0,j}=f_{0}^{\otimes j}.
\Eq(2.26)
$$
Such a notion of statistical independence, which corresponds to a complete
factorization of the RDM's, is not compatible (but for the condensed Bose
state) with the B-E and F-D statistics which exhibit intrinsic
correlations even for non interacting particle systems.
States describing free Bosons or Fermions are usually called quasi-free
and are defined in terms of the RDM's by the following formula:
$$
\r_{j}(x_{1},\dots,x_{j};y_{1},\dots,y_{j})
=
\sum_{\pi\in\Cal P_{j}}\th^{s(\pi)} \prod_{i=1}^{j}\r(x_{i}, y_{\pi(i)}),
\Eq(qf)
$$
for some positive definite operator $\r$ on
$L_{2}(\Bbb R^{3})$ with kernel $\r(x,y)$.
We show in Appendix how to construct explicitly quasi-free states for Bosons.

From now on we assume that the initial sequence \equ(is) for the 
rescaled problem
\equ(rp) is given by the Wigner transform of a quasi-free state \equ(qf)
generated by $\r(x,y)=\r^{\e}(x,y)$ given by \equ(r). As a consequence the
initial sequence $\{f_{j}^{0}\}_{j=1}^{\infty}$ for
the hierarchy \equ(rp)
is of the form
$$
f_{j}^{0}(X_{j},V_{j})
=\sum_{\pi\in\Cal P_{j}}\th^{s(\pi)}f_{j}^{\pi}(X_{j},V_{j}),
\Eq(qf1)
$$
with
$$
\eqalign{&
f_{j}^{\pi}(x_{1},\dots,x_{j},v_{1},\dots,v_{j})
={1\over (2\pi)^{3j}}
\int dy_{1}\dots\int dy_{j}\int dw_{1}\dots\int dw_{j}
\cr&
\prod_{k=1}^{j}
\left(
\ex^{iy_{k}\cdot v_{k}
+i w_{k}\cdot {x_{k}-x_{\pi(k)} \over \e}
-i w_{k}\cdot {y_{k}+y_{\pi(k)} \over 2}}
f_{0}\Big({x_{k}+x_{\pi(k)}\over 2}-{\e\over
4}(y_{k}-y_{\pi(k)}),w_{k}\Big)\right).     
}
\Eq(pi)
$$
Eq. \equ(pi) follows from \equ(qf) and \equ(r).

We underline once more that, in the present approach, the dynamics is
given by the hierarchy of equations \equ(rp) which are completely
equivalent to the Schr\"odinger equation, while the statistics enters only
in the structure of the initial state.

In the weak-coupling limit $\e \to 0$, we expect that $f^\e_j(t)$
converges to a factorized state (because the effects of statistics
disappear in the macroscopic limit). On the more  each factor 
should be 
solution to the U-U
equation  (the collisions being affected by the statistics because they
involve microscopic scales).

\bigskip\bigskip


\heading 3. The main result.\endheading

\numsec=3
\reset
\numfor= 1
\bigskip


Let $f=f(x,v,t)$ be a solution to the U-U equation and set
$f_j(\,\cdot\,,\,\cdot\,,t)= f^{\otimes j}(\,\cdot\,,\,\cdot\,,t)$. Then
the sequence $\{f_j\}_{j=1}^\infty$ satisfies the following hierachy of
equations:
$$
(\pt_t+\sum_{i=1}^j v_i\cdot \n_{x_i})f_j
= Q_{j,j+1} f_{j+1} + Q_{j,j+2}f_{j+2}.
\Eq(3.1)
$$
Here the $Q_{j,j+1}$ contribution, a "two particles term"
in the terminology used below, is
$$
\eqalign{
(Q_{j,j+1} f_{j+1})&(X_j,V_j)
=
\sum_{k=1}^j\int dv_k'\int dv_{j+1}\int dv_{j+1}'\,
W(v_k,v_{j+1}|v_k',v_{j+1}')\cr&
\Big\{f_{j+1}(X_j, x_k;v_1,\dots, v_k',\dots v_{j+1}')
-
f_{j+1}(X_j, x_k;v_1,\dots,  v_{j+1})\Big\},}
\Eq(3.2)
$$
and the $Q_{j,j+2}$ contribution, a "three particles term", is
$$
\eqalign{&
(Q_{j,j+2} f_{j+2})(X_j,V_j)=
8\pi^3\th \sum_{k=1}^j\int 
dv_k'\int dv_{j+1} \int
dv_{j+1}'W(v_k,v_{j+1}|v_k',v_{j+1}') \, \cr
\Bigg\{
&f_{j+2}(X_j, x_k, x_k;v_1,\dots, v_k',\dots v_{j+1}',v_k)
+
f_{j+2}(X_j, x_k, x_k;v_1,\dots, v_k',\dots v_{j+1}',v_{j+1})
\cr
&
\qquad
-
f_{j+2}(X_j, x_k\, x_k;v_1,\dots,  v_{j+1},v_k')
-
f_{j+2}(X_j\, x_k, x_k;v_1,\dots,  v_{j+1},v_{j+1}')
\Bigg\}.}
\Eq(3.3)
$$
Also, $(X_n, y)$ denotes the $(n+1)$-sequence $(x_1,\dots,x_n,y)$.

A formal solution to the hierarchy \equ(3.1) is given by the
following series expansion:
$$
\eqalign{
f_j(t)=&\sum_{n\ge 0}\sum_{ \a_1\dots\a_n \atop \a_i=1,2}
\int_0^tdt_1\int_0^{t_1}dt_2\dots\int_0^{t_{n-1}}dt_n
\cr&
S(t-t_1)Q_{j,j+\a_1}S(t_1-t_2)Q_{j+\a_1,j+\a_1+\a_2}S(t_2-t_3)\dots
\cr&
Q_{j+\a_1+\dots+\a_{n-1}, j+\a_1+\dots+\a_{n}}S(t_n)
f_0^{\otimes(j+\a_1+\dots+\a_{n})},}
\Eq(3.4)
$$
where $S(t)$ denotes the fream stream operator, namely,
$$
(S(t)f_j)(X_j,V_j)= f_j(X_j-V_jt, V_j).
\Eq(3.5)
$$

As for the solution
to the $\e$-dependent  hierarchy \equ(rp),
we can also expand $f^j_\e$ in the similar way, 
at least at the formal level. This gives
$$
\eqalign{
f^\e_j(t)=&\sum_{n\ge 0}\sum_{ \g_1\dots\g_n \atop \g_i=0,1}
\int_0^tdt_1\int_0^{t_1}dt_2\dots\int_0^{t_{n-1}}dt_n
\cr&
S(t-t_1)P^\e_{j,j+\g_1}S(t_1-t_2)P^\e_{j+\g_1,j+\g_1+\g_2}S(t_2-t_3)\dots
\cr&
P^\e_{j+\g_1+\dots+\g_{n-1}, j+\g_1+\dots+\g_{n}}S(t_n)
f^0_{j+\g_1+\dots+\g_{n}},}
\Eq(3.6)
$$
where $f^0_j$ is an initial quasi-free state given by \equ(pi), and
$$
P^\e_{j,j+1}=\e^{-{7\over 2}} C_{j+1}^\e, \quad
P^\e_{j,j}=\e^{-{1\over 2}}T^\e_j.
\Eq(3.7)
$$

We are not able to show the convergence of $f^\e_j(t)$ to $f_j(t)$ in the
limit $\e\to 0$ even for short times. However we are going to show that
the two series agree up to the second order in the potential. Namely,
defining the second order contributions
$$
g(t):=S(t)f_0+ \int_0^tdt_1 \; S(t-t_1)Q_{1,2}S(t_1)f_0^{\otimes 2} 
+
\int_0^tdt_1 \; S(t-t_1)
Q_{1,3}S(t_1)f_0^{\otimes 3},
\Eq(3.8)
$$
associated with $f_j(t)$, and
$$
\eqalign{g^\e(t):=S(t)f_0&
+\e^{-{7\over 2}}\int_0^tdt_1 \; S(t-t_1)C^\e_{2}S(t_1) f^0_2
\cr&
+ 
\e^{-4}\int_0^tdt_1 
\; \int_0^{t_1}dt_2 \; S(t-t_1)C^\e_{2}S(t_1-t_2)T^\e_2S(t_2)f^0_2
\cr&
+
\e^{-7}\int_0^tdt_1\; 
\int_0^{t_1}dt_2\; S(t-t_1)C^\e_{2}S(t_1-t_2)C^\e_3S(t_2)f^0_3,}
\Eq(3.9)
$$
associated with $f_j^\e(t)$, 
we rigorously prove below the convergence of $g^\e(t)$ to 
$g(t)$, under suitable assumptions
on the data of the problem.

\medskip

\noi
$\underline{\text{Assumptions}}$:
We require $\phi$ to be real and even, so that  $\widehat
\phi$ is real. In particular 
$$
\widehat \phi(k)= \overline{\widehat\phi(-k)}=\widehat\phi(-k).
$$ 
This the most important assumption we need in the analysis.
Besides, we shall need to deal with "smooth" data. 
Quantitatively, we assume the following regularity:
$$
(1+|\xi| )^\a \sum_{|\b|\le 2}|D_\x^\b\widehat \phi(\xi)|\in L_1,
$$
for a sufficiently large $\a$, and
$$
\eqalign{& 
f_0(x,v)\in L_1,\cr&
(1+|\xi|+|\eta|)^\a \sum_{0\le|\b|\le 2}
\sum_{0\le|\g|\le 2}|(D_\x^\b + D_\h^\g)\widehat f_0(\xi,\eta)|\in L_1,}
\Eq(3.10)
$$ 
for a sufficiently large $\a$ as well.
In \equ(3.10), $\b$ and $\g$ denote 
multi-indices, and $D^\b_\x$, $D^\g_\h$ denote derivatives with respect 
to the variables $\x$ and $\h$.
Note that throughout this paper
we use the following normalization for the Fourier transform:
$$
\eqalign{&
\widehat f(h)=(\Cal F_xf)(h)=\int_{\Bbb R^n} dx \; f(x) \, \ex^{-ih\cdot x},
\cr&
f(x)=
(\Cal F^{-1}_h \widehat f)(h)=
{1\over (2\pi)^n}\int _{\Bbb R^n}d h  \; \widehat f(h) \,
\ex^{ih\cdot x}.}
\Eq(2.22)
$$
\hfill\qed

\bigskip

\noi
Our main result is the

\proclaim{Theorem}
Under the above assumptions, we have
$$
\lim_{\e\to 0} \widehat g^\e(\x,\h,t)=\widehat g(\x,\h,t),
\quad \text{for any $t>0$ and any $(\xi,\eta)\in \Bbb R^6$}.
$$
\endproclaim

\bigskip

\noindent
$\underline{\text{Remark}}$:
In the above statement (and the proofs given below), we found 
convenient 
to treat the terms in \equ(3.8) 
and \equ(3.9) in terms of their Fourier transforms,
for which the convergence arises more naturally.
However, we would like to stress that in the companion paper 
[\rcite{BCEP}],
a stronger, but analogous, result is formulated in terms of
the pointwise convergence in the $x-v$ space, hence without going
to the Fourier space.
\hfill\qed

\bigskip

Before entering the details of the proof we first analyze all the 
contributions in the 
right hand side of \equ(3.9).

The two-particle terms are
(we skip the unessential operator
$\int_0^tdt_1 \; S(t-t_1)$)
$$
\Cal S^\pi_{2,\e}:= \e^{-{7\over 2}}C_2^\e S(t_1)f_2^\pi,
\Eq(3.11)
$$
where the permutation $\pi$ may take the two values $\pi=(1,2)$ 
or $\pi=(2,1)$, together with
$$
\Cal T_{2,\e}^\pi:= \e^{-4} \int_0^{t_1}dt_2 \; C_2^\e S(t_1-t_2)T_2^\e 
S(t_2)f_2^\pi,
\Eq(3.12)
$$
with $\pi$ taking the values $\pi=(1,2)$ or $\pi=(2,1)$.
There are four such terms.

The three-particle terms are twelve, namely:
$$
\Cal W^\pi_{3,\e}:=\e^{-7} \int_0^{t_1}dt_2 \; C_{1,2}^\e 
S(t_1-t_2)C_{1,3}^\e
S(t_2)f_3^\pi,
\Eq(3.13)
$$
and  
$$
\Cal V^\pi_{3,\e}:=\e^{-7} \int_0^{t_1}dt_2 \; C_{1,2}^\e 
S(t_1-t_2)C_{2,3}^\e S(t_2)f_3^\pi,
\Eq(3.14)
$$
with $\pi\in \Cal P_3$, the set of the permutations of three objects, whose
cardinality is six.

Note that all the above terms are funtions of $(x,v)$ (and $t_1$ of 
course).

For further convenience, and in view of the proof of our main result,
we readily express all these contributions in terms of their Fourier
transforms. 

We start with the following obvious three formulae for the basic
operators $S(t)$, $T_2$, and $C_2$ (see \equ(3.5), \equ(2.18)-\equ(2.19),
and \equ(2.20)-\equ(2.21), respectively):
$$
\eqalign{
&
\widehat T_2^\e \widehat f(\xi_1,\xi_2;\eta_1,\eta_2)
=-i\sum_{\s=\pm 1}\s\int dh \;
{\widehat\phi(h)\over (2\pi)^3} \;
\ex^{i\s{h\over 2}\cdot (\eta_2-\eta_1)}
\widehat f\left(\xi_1-{h\over \e}, \xi_2+{h\over \e};\eta_1,\eta_2\right),
\cr&
\widehat C_2^\e \widehat f(\xi;\eta)
=\widehat T_2^\e \widehat f(\xi,0;\eta,0)
= -i\sum_{\s=\pm 1}\s\int dh \;
{\widehat\phi(h)\over (2\pi)^3} \; 
\ex^{-i\s{h\over 2}\cdot \eta}
\widehat f\left(\xi-{h\over \e}, {h\over \e};\eta,0\right),
\cr&
\widehat S(t) \widehat f(\xi;\eta)
=\widehat f(\xi,\eta+\xi t).
}
$$
These relations give in \equ(3.11) through \equ(3.14):
$$
\widehat{\Cal S}^\pi_{2,\e}(\xi,\eta)=
-i{\e^{-{7\over 2}}\over (2\pi)^3}
\sum_{\s=\pm 1}\s\int dh \; 
\widehat\phi(h) \;
\ex^{-i{\s\over 2}h\cdot\h } \;
\widehat f_2^\pi\left(\x-{h\over \e}, {h\over \e};
\h+t_1(\x-{h\over \e}), t_1{h\over \e}\right).
\Eq(3.17')
$$
$$
\eqalign{
\widehat{\Cal T}^\pi_{2,\e}(\xi,\eta)
=&
-{\e^{-4}\over (2\pi)^6}
\sum_{\s_1,\s_2=\pm 1}\s_1\s_2
\int_0^{t_1}dt_2 \int dh_1\int dh_2 \;
\widehat\phi(h_1)\, \widehat\phi(h_2)\, \cr&
\ex^{-i{\s_1\over 2}h_1\cdot \eta} \,
\ex^{-i{\s_2\over 2}h_2\cdot (\eta+\xi(t_1-t_2)-{2\over \e}(t_1-t_2)h_1)}
\cr&
\widehat f_2^\pi\left(\xi-{1\over \e}(h_1+h_2), {1\over \e} 
\left(h_1+h_2\right);
\eta+t_1
\xi-{t_1h_1+t_2h_2\over \e}, {t_1h_1+t_2h_2\over
\e}\right),}
\Eq(3.18)
$$
$$
\eqalign{
\widehat{\Cal W}^\pi_{3,\e}(\xi,\eta)=&
-{\e^{-7}\over (2\pi)^6}
\sum_{\s_1,\s_2=\pm 1}\s_1\s_2
\int_0^{t_1}dt_2 \int dh_1\int dh_2 \;
\widehat\phi(h_1) \, \widehat\phi(h_2)
\cr&
\ex^{-i{\s_1\over 2}h_1\cdot \eta} \, 
\ex^{-i{\s_2\over 2}h_2\cdot (\eta+(\xi-{h_1\over \e}) (t_1-t_2))}
\cr&
\widehat f_3^\pi\left(\xi-{1\over \e}(h_1+h_2), {h_1\over \e},{h_2\over \e}; 
\eta+t_1\xi-{t_1h_1+t_2h_2\over \e}, {t_1h_1\over \e}, {t_2 h_2\over \e},
\right),}
\Eq(3.19)
$$
$$
\eqalign{
\widehat{\Cal V}^\pi_{3,\e}(\xi,\eta)=&
-{\e^{-7}\over (2\pi)^6}
\sum_{\s_1,\s_2=\pm 1}\s_1\s_2
\int_0^{t_1}dt_2 \int dh_1\int dh_2 \;
\widehat\phi(h_1) \, \widehat\phi(h_2)
\cr&
\ex^{-i{\s_1\over 2}h_1\cdot \eta} \, 
\ex^{-i{\s_2\over 2}h_2\cdot {h_1\over \e}(t_1-t_2)}
\cr&
\widehat f_3^\pi\left(\xi-{h_1\over \e}, {h_1-h_2\over \e},{h_2\over \e}; 
\eta+t_1\xi-{t_1h_1\over \e}, {t_1h_1-t_2h_2\over \e}, {t_2 h_2\over \e}
\right).}
\Eq(3.20)
$$

Starting form those expressions,
the plan of the proof is the following. In Section 4 we evaluate the two
particle terms $\Cal S^\pi_{2,\e}$ and $\Cal T^\pi_{2,\e}$. We prove
that they converge towards the associated two particles terms in the U-U
equation. In Section 5 we deal with the three-particle terms associated
to the permutations $\pi$ with a fixed element. Those are shown to
converge towards the associated three particles terms in the U-U equation,
while contributing by the quantity $\widehat \phi(v'-v)^2+\widehat
\phi(v'-v_*)^2$ to the transition kernel $W$ (see \equ(W)). Finally in
Section 6 we treat the three particle terms relative to cyclic
permutations. We recover in this way the missing contribution to the
transition kernel, namely the cross term $\th \widehat \phi(v'-v) \,
\widehat \phi(v'-v_*)$.

For sake of simplicity we shall carry out the computations for Bosons
($\th =1$), being clear that the Fermionic case is just the same with
suitable changes of sign.

\bigskip\bigskip


\heading 4. Two-particle terms.\endheading

\numsec= 4
\reset
\numfor= 1
\bigskip


\noi
We introduce the partial Fourier transform 
$$
\widetilde f_j^\pi(x_1,\dots,x_j;\eta_1,\dots,\eta_j)
:= \int dv_1\dots\int dv_j \;
\ex^{-i\sum_{k=1}^j v_k\cdot \h_k} \,
f_j^\pi(x_1,\dots,x_j;v_1,\dots,v_j).
\Eq(fv)
$$
As a consequence of \equ(pi) we have
$$
\eqalign{
\widetilde f_j^\pi(x_1,\dots,x_j&;\eta_1,\dots,\eta_j)
=
\cr&
\prod_{k=1}^j
\widetilde f_0\left(
{x_k+x_{\pi(k)}\over 2}  -  {\e}{\h_k-\h_{\pi(k)}\over 4}, 
-   {x_k-x_{\pi(k)}\over \e}+{\h_k+\h_{\pi(k)}\over 2}\right).
}
\Eq(til)
$$
In particular, we have the obvious 
$$
\widetilde f_2^{(1,2)}(x_1,x_2;\eta_1,\eta_2)
=\widetilde f_0(x_1,\eta_1)\widetilde f_0 (x_2,\eta_2),
\Eq(4.1)
$$
together with
$$
\eqalign{
\widetilde f_2^{(2,1)}(x_1,x_2;\eta_1,\eta_2)=&
\widetilde f_0
\left(
{x_1+x_2\over 2} -  {\e\over 4}(\eta_1-\eta_2);
-  {x_1-x_2\over\e}+{\eta_1+\eta_2\over 2}
\right)\cr&
\widetilde f_0
\left(
{x_1+x_2\over 2} +  {\e\over 4}(\eta_1-\eta_2);
 +  {x_1-x_2\over \e}+{\eta_1+\eta_2\over 2}
\right).}
\Eq(4.2)
$$
Hence, upon now performing the complete Fourier transform,
we obtain,
$$
\widehat f_2^{(1,2)}(\xi_1,\xi_2;\eta_1,\eta_2)=
\widehat f_0(\x_1,\eta_1)
\widehat  f_0(\x_2,\eta_2),
\Eq(4.3)
$$
together with
$$
\eqalign{\widehat f_2^{(2,1)}(\x_1,\x_2;\eta_1,\eta_2)=&
\e^3\int dy_1\int dy_2 \;
\ex^{-i\x_1\cdot(y_1+{\e\over 2}y_2)} \,
\ex^{-i\x_2\cdot(y_1-{\e\over 2}y_2)} \,
\cr&
\widetilde f_0
\left(
y_1 -  {\e\over 4}(\eta_1-\eta_2);
-  y_2+{\eta_1+\eta_2\over 2}
\right)
\cr&
\widetilde f_0
\left(
y_1 +  {\e\over 4}(\eta_1-\eta_2);
 + y_2+{\eta_1+\eta_2\over 2}
\right).}
\Eq(4.4)
$$

\bigskip

We are now in position to analyse the term ${\Cal S}_{2,\e}^{\pi}$
for $\pi=(1,2)$ and $\pi=(2,1)$.
First, 
using the identity
$$
\sum_{\s=\pm 1}\s \ex^{-i\s {h\over 2}\cdot \eta}
=-2 i\sin{h\cdot\eta\over 2},
$$
we get the 
the explicit expression:
$$
\widehat{\Cal S}_{2,\e}^{\pi}=
-{2\e^{-{7\over 2}}\over (2\pi)^3}\int dh \, \widehat\phi(h)
\sin\left({h\cdot\eta\over 2}\right)
\widehat f_2^{\pi}
\left(
\xi-{h\over \e},{h\over\e};
\eta+(\x-{h\over \e})t_1,{h\over \e}t_1\right).
\Eq(4.5)
$$
In the case 
of $\widehat{\Cal S}_{2,\e}^{(1,2)}$,
a change of
variable $h \to \e h$ then gives, using \equ(4.3), the value
$$
\widehat{\Cal S}_{2,\e}^{(1,2)}(\x,\eta)=
-{2\e^{-{1\over 2}}\over (2\pi)^3}\int dh \,
\widehat\phi(\e h)\sin\left(\e{h\cdot\eta\over 2}\right)
\widehat f_0(\xi-{h};\eta+(\x-h)t_1)
\widehat f_0({h},{h}t_1).
\Eq(4.6)
$$
Therefore, we may estimate
$$
\eqalign{&|\widehat{\Cal S}_{2,\e}^{(1,2)}| 
\le \cr 
&\le C\sqrt\e \, 
\|\widehat\phi\|_{L_\infty}
\int dh\, |h|\,\left( |\eta + (\xi -h)t_1| +
|\xi - h|t_1\right) |\widehat f_0(h;ht_1)|\,
|\widehat f_0|( \xi - h ; \eta + (\xi -h)t_1)\cr
& \le C\sqrt\e \left(
\sup_{\xi,\eta} |\eta|\/ |\widehat f_0|(\xi;\eta) 
\int d\xi \sup_{\eta} |\xi|\/ |\widehat f_0|(\xi;\eta) +t_1
\sup_{\xi,\eta}  |\xi| \/|\widehat f_0|(\xi;\eta) 
\int d\xi |\x|\sup_{\eta}  \/|\widehat f_0(\xi;\eta)| \right)
\cr
&\le C
\sqrt\e
\sup_{\xi,\eta} 
\left(
(|\xi| + |\eta|)\,|\widehat f_0|(\xi;\eta)
\right)
\int d\xi\, \sup_{\eta}
\left(
(|\xi| + |\eta|)\,|\widehat f_0|(\xi;\eta)
\right),
}
\Eq(4.7)
$$ 
and the corresponding contribution vanishes with $\e$.

In the case of $\widehat{\Cal S}_{2,\e}^{(2,1)}$
on the other hand, 
equations 
\equ(4.5) and \equ(4.4) give 
$$
\eqalign{
\widehat{\Cal S}_{2,\e}^{(2,1)}(\x,\eta)
&=-{2\e^{-{1\over 2}}\over (2\pi)^3}
\int dh\int dy_1 \int dy_2 \; \widehat\phi(h) \,
\sin\left({h\cdot\eta\over 2}\right)
\ex^{-i(y_1+{\e\over 2}y_2)\cdot \x} \, \ex^{ih\cdot y_2}
\cr&
\quad \widetilde f_0\left(y_1-{\e\over 4}\left[\eta+\x t_1-{2h\over 
\e t_1}\right];
-y_2+{\h+\x t_1\over 2}\right)
\cr&
\quad \widetilde f_0\left(y_1+{\e\over 4}\left[\eta+\x t_1-{2h\over 
\e t_1}\right];
y_2+{\h+\x t_1\over 2}\right) \,
\cr&
= - {2\e^{-{1\over 2}}\over (2\pi)^3}
\int dh\int dy_1 \int dy_2 \;
\widehat\phi(h) \, \sin\left({h\cdot\eta\over 2}\right) \,
\ex^{-i(y_1+\e{y_2\over 2})\cdot\x}\ex^{ih\cdot y_2}
\cr&
\qquad
\widetilde f_0\left(y_1+{h\over 2}t_1;-y_2+{\h+\x t_1\over 2}\right)
\widetilde f_0\left(y_1-{h\over 2}t_1;y_2+{\h+\x t_1\over 2}\right)
+O(\sqrt\e).
}
\Eq(4.8)
$$
By the parity of $\phi$,
the first term in the right hand side is vanishing:
Indeed, it is
antisymmetric in the exchange $h\to -h$ and $y_2\to -y_2$.
Note that the mechanism that makes the dominant, $O(\e^{-1/2})$, 
contribution of $\widehat{\Cal S}_{2,\e}^{(2,1)}$, vanish in
the limit, is very different 
from the one involved in the vanishing of $\widehat{\Cal S}_{2,\e}^{(1,2)}$:
here, antisymmetry plays a crucial role. This aspect will play an even 
more important, and more intricate, role in the next two sections.
There remains to prove that the 
$O(\sqrt \e)$ term in \equ(4.8) indeed has the claimed size. It
can be written as
$$
\eqalign{&
- {2\e^{{1\over 2}}\over (2\pi)^6}
\int dh \, dy_2 \, d\xi_1 \
{\widehat \phi}(h) \,\sin\left({h\cdot\eta\over 2}\right) \,
\widehat f_0\left( \xi_1;  {\eta + \xi t_1 \over 2 } - y_2\right)\,
\widehat f_0\left( \xi - \xi_1;  {\eta + \xi t_1 \over 2 } + y_2\right)
\cr&
\qquad
\ex^{ih\cdot y_2 + i {t_1 \over 2} h \cdot (2\xi_1 -\xi) }
\left( 1 - \ex^{i{\e \over 2} \left( -y_2 \cdot \xi +(\xi-2 \xi_1) \cdot 
\left( {\eta + \xi t_1 \over 2 } \right) \right)}
\right).
}
$$
It may be estimated by
$$
\eqalign{&
C\e^{{1\over 2}} \int dh\,{|\widehat \phi|}(h)
\int dy_2 \,d\xi_1 \,
|\widehat f_0|\left( \xi_1;  {\eta + \xi t_1 \over 2 } + y_2\right)\,
|\widehat f_0|\left( \xi - \xi_1;  {\eta + \xi t_1 \over 2 } - y_2\right)
\cr&
\left(|\xi_1|\,\left|{\eta + \xi t_1 \over 2 } + y_2\right| + 
|\xi - \xi_1|\, \left|  {\eta + \xi t_1 \over 2 } - y_2\right|\right).
}
$$
Therefore the term $\Cal S_{2,\e}^{(2,1)}$ vanishes as well.

As a conclusion, all terms $\Cal S_{2,\e}^\pi$, which
are the ones that are linear in $\phi$, vanish
in the limit $\e\to 0$.

\bigskip

We now pass to the evaluation of the terms $\Cal T_{2,\e}^\pi$. 

The contribution
$\Cal T_{2,\e}^{(1,2)}$ has been already considered in Ref.
[\rcite{BCEP}]. However, for sake of completeness,
we analyze this term in the present context as
well.
Using \equ(4.3) in \equ(3.18), and performing the change of variables:
$$
{h_1+h_2\over \e}=k, \quad 
h_2=h, \quad
{t_1-t_2\over \e}=s,
\Eq(4.9)
$$
we arrive at
$$
\eqalign{\widehat{\Cal T}_{2,\e}^{(1,2)}(\x,\eta)&=
-{1\over (2\pi)^6}
\sum_{\s_1,\s_2=\pm 1}\s_1\s_2
\int_0^{t_1\over \e}ds\int dh\int dk \; 
\widehat\phi(h)\widehat\phi(-h+\e k)
\cr&
\ex^{-i[\eta\cdot h\left({\s_2-\s_1\over 2}\right)+\s_2h^2s]} \, 
\ex^{-i\e{\s_1\over 2}k\cdot \eta}\, 
\ex^{-i\e\s_2{h\over 2}\cdot [\x s-2sk]}
\cr&
\widehat f_0(\x-k;\eta+t_1\x+hs-kt_1) \,
\widehat f_0(k;kt_1-hs).}
\Eq(4.10)
$$
This term converges formally to
$$
\eqalign{
\widehat{\Cal T}_{2}^{(1,2)}(\x,\eta)&=
-{1\over (2\pi)^6}\sum_{\s_1,\s_2=\pm 1}\s_1\s_2
\int_0^{+\infty}ds\int dh\int dk \; |\widehat\phi(h)|^2
\cr&
\ex^{-i[\eta\cdot h\left({\s_2-\s_1\over 2}\right)+\s_2h^2s]} \,
\widehat f_0(\x-k;\eta+t_1\x+hs-kt_1) \,
\widehat f_0(k;kt_1-hs).}
\Eq(4.11)
$$
To justify the limit we split the integration in $ds$ over the two
intervals $[0,1]$, $[1,+\infty]$. In the first interval we bound the
integrand by
$$
\|\widehat\phi\|_{L_\infty} \,
\|\widehat f_0\|_{L_\infty} \,
|\widehat\phi(h)|\sup_\eta \,
|\widehat f_0(k,\eta)|,
\Eq(4.12)
$$ 
which is a $L_1(dk\,dh)$ function for any
$s\in [0,1]$. In the second part of the integration domain, after the
change of variables $h\to (kt_1-hs)$, we bound the integrand by
$$
\|\widehat\phi\|^2_{L_\infty} \,
\|\widehat f_0\|_{L_\infty} \,
|\widehat f_0(k,\eta)|{1\over s^3},
\Eq(4.12')
$$ 
which is a $L_1(dk\,dh\,ds)$ function on
$\Bbb R^3\times\Bbb R^3\times (1,+\infty)$.
The claimed convergence in \equ(4.11)
is then 
consequence of the 
Dominated Convergence Theorem. It holds uniformly in $\x$, $\eta$.

We now evaluate $\widehat{\Cal T}_{2,\e}^{(2,1)}$.
Inserting \equ(4.4) 
in \equ(3.18), and rescaling time
$t_1-t_2=\e s$, the resulting expression is:
$$
\eqalign{
\widehat{\Cal T}_{2,\e}^{(2,1)}(\x,\eta)=
&
-{1\over (2\pi)^6}\sum_{\s_1,\s_2=\pm 1} \s_1\s_2
\int_0^{t_1\over \e} ds \int dh_1 \int dh_2\int dy_1\int dy_2 \;
\widehat\phi(h_1) \, \widehat\phi(h_2)
\cr&
\ex^{-i {\s_1\over 2} h_1\cdot\eta} \,
\ex^{-i {\s_2\over 2} h_2\cdot(\eta+\e s\x- 2 s h_1)} \,
\ex^{-i \x\cdot(y_1+{\e\over 2} y_2)} \,
\ex^{i(h_1+h_2)\cdot y_2} \,
\cr&
\widetilde f_0
\left(y_1-{\e\over 4}\left[\eta+\x t_1-2{h_1t_1+h_2(t_1-\e s)\over 
\e}\right]; 
-y_2+{\eta+\x t_1\over 2}\right)
\cr&
\widetilde f_0
\left(y_1+{\e\over 4}\left[\eta+\x t_1-2{h_1t_1+h_2(t_1-\e s)\over 
\e}\right];
y_2+{\eta+ \x t_1\over 2}\right).
}\Eq(4.13)
$$
Now the formal limit is:
$$
\eqalign{
\widehat{\Cal T}_{2}^{(2,1)}(\x,\eta)=
&
-{1\over (2\pi)^6}\sum_{\s_1,\s_2=\pm 1} \s_1\s_2\int_0^{+\infty} ds
\int dh_1 \int dh_2\int dy_1\int dy_2 \;
\widehat\phi(h_1) \, \widehat\phi(h_2)
\cr&
\ex^{-i {\s_1\over 2} h_1\cdot\eta} \,
\ex^{-i {\s_2\over 2} h_2\cdot(\eta- 2 s h_1)} \,
\ex^{-i \x\cdot y_1} \,
\ex^{i(h_1+h_2)\cdot y_2} \,
\cr&
\widetilde f_0
\left(y_1+t_1{h_1+h_2\over 2}; -y_2+{\eta+\x t_1\over 2}\right) \,
\widetilde f_0
\left(y_1-t_1{h_1+h_2\over 2}; y_2+{\eta+\x t_1\over 2}\right).
}\Eq(4.14)
$$
To justify the limit we have to show the uniform (with respect to $\e$) 
integrability  of the integrand in the right hand side of \equ(4.13). To
outline the decay with respect to the $s$ variable we observe
$$
\ex^{i\s_2h_1\cdot h_2s}= -{1\over s^2|h_2|^4}(h_2\cdot\n_{h_1})^2 
\ex^{i\s_2 h_1\cdot h_2 s},
\Eq(4.15)
$$
and then proceed with the natural integration by parts with
respect to $h_1$ in \equ(4.13) (Recall that $1/|h_2|^2$ is integrable
close to the origin in dimension $d=3$).
Splitting the integral in $ds$ as before, 
we may apply the
Dominated Convergence Theorem,
upon using the smoothness of $\widehat \phi$ 
and $\widetilde f_0$, thus justifying the above formal limit.

\bigskip

Our last task is to interpret the result we have obtained,
in terms of the U-U equation. To do so, we first go back
to the original variables, expressing $\Cal T_{2}^{(1,2)}$ and
$\Cal T_{2}^{(2,1)}$ as functions of $(x,v)$. A 
straightforward computation yields, on the one hand,
$$
\eqalign{
{\Cal T_{2}^{(1,2)}}(x,v)&=
-{1\over(2\pi)^3}\sum_{\s_1,\s_2=\pm 1}\s_1\s_2
\int dh\int dv_1\int dv_2 \;
\delta\left(v-v_1-h{\s_2-\s_1\over 2}\right) \,
\cr&
\Delta^+\left(-h\cdot(\s_2h +  (v_1-v_2)\right)\,
|\widehat\phi(h)|^2 \,
f_0(x-v_1t_1,v_1) \, f_0(x-v_2t_1,v_2),
}
\Eq(4.16)
$$
and, on the other hand (with $h_2 = h$),
$$
\eqalign{
{\Cal T_{2}^{(2,1)}}(x,v)&=
-{1\over(2\pi)^3} \sum_{\s_1,\s_2=\pm 1}\s_1\s_2
\int dh\int dv_1\int dv_2
\cr&
\delta\left( v-v_1{1-\s_1\over 2}
-v_2{1+\s_1\over  2}-h{\s_2-\s_1\over2}
\right)
\cr&
\Delta^+(\s_2 h\cdot(-h - v_1+ v_2 ))\,
\widehat\phi(h) \, \widehat\phi( -  h-v_1+v_2) \,
\cr&
f_0(x-v_1t_1,v_1) \, f_0(x-v_2t_1,v_2).
}
\Eq(4.17)
$$
Here we define the distribution
$$
\Delta^+(\a):=
\int_0^\infty ds \; \ex^{i\a s}
=\pi\d(\a) +i \text{\rm p.v.} \left({1\over \a}\right).
\Eq(4.18)
$$
Now, on both preceding formulae, we readily observe the following 
important fact. The parity of 
$\widehat\phi$, and the symmetries $h\to -h$, $\s_1\to -\s_1$,
$\s_2\to -\s_2$ in
\equ(4.16), and $h\to -h$, $\s_1\to -\s_1$, $\s_2\to -\s_2$, 
$v_1\leftrightarrow v_2$ in
\equ(4.17), show that $\Delta^+$ may be replaced by $\pi\d$ everywhere. 
This will eventually give, as shown next, the desired conservation of
energy in the limiting U-U equation.

There remains to actually perform the sum $\sum \s_1 \s_2$ in 
\equ(4.16) and \equ(4.17), in order to identify the very value
of $\Cal T_{2}^{(1,2)}$ and $\Cal T_{2}^{(2,1)}$.
For $\Cal T_2^{(1,2)}$ we make the following 
choice:
\bigskip
\centerline{
\vbox{\tabskip=0pt \offinterlineskip
\def\tablerule{\noalign{\hrule}}
\halign to250pt{
\strut#& \vrule#\tabskip=1em plus2em&
\hfil#&\vrule#& \hfil#\hfil& \vrule#&
\hfil#&\vrule#& \hfil#\hfil& \vrule#&
\hfil#& \vrule#
\tabskip=0pt
\cr
\tablerule
&&$\s_1$&&$\s_2$&&$h\hskip.4cm$&&$v_1$&&$v_2$&\cr
\tablerule 
&&$1$&&$1$&&$v'-v$&&$v$&&$v_*$&\cr
\tablerule 
&&$-1$&&$-1\hskip.27cm$&&$v-v'$&&$v$&&$v_*$&\cr
\tablerule 
&&$1$&&$-1\hskip.27cm$&&$ v'-v $\hskip-.2cm&&$v'$&&$v_*'$&\cr
\tablerule 
&&$-1$&&$1$&&$ v-v' \hskip-.2cm$&&$v'$&&$v_*'$&\cr
\tablerule 
}} 
}
\bigskip
\noi
This results in the final expression:
$$
\eqalign{
{\Cal T_{2}^{(1,2)}}(x,v)=&
{1\over 4\pi^2}\int dv_*\int dv'\int dv_*' \;
\delta(v_*+v-v_*'-v')
\cr&
\delta\left({1\over2}(v_*\/^2+v^2-v_*'\/^2-v'\/^2)\right) \,
|\widehat\phi(v'-v)|^2 \,
(f'f'_*-ff_*),
}
\Eq(4.19)
$$
where, with abuse of notation, we set the "transported quantities"
$$
f=f_0(x-vt_1,v), \quad
f_*=f_0(x-v_*t_1,v_*), \quad 
f'=f_0(x-v't_1,v'), \quad
f_*'=f_0(x-v_*'t_1,v_*').
$$
Notice that,
by changing $v'\leftrightarrow v'_* $,
and using the conservation 
of momentum, we
may replace 
${1\over 2}|\widehat\phi(v'-v)|^2$
by ${1\over 2}|\widehat\phi(v'-v_*)|^2$ in \equ(4.19).

\noi
Besides, for $\Cal T_2^{(2,1)}$ we make the following changes
of variables:
\bigskip
\centerline{
\vbox{\tabskip=0pt \offinterlineskip
\def\tablerule{\noalign{\hrule}}
\halign to250pt{
\strut#& \vrule#\tabskip=1em plus2em&
\hfil#&\vrule#& \hfil#\hfil& \vrule#&
\hfil#&\vrule#& \hfil#\hfil& \vrule#&
\hfil#& \vrule#
\tabskip=0pt
\cr
\tablerule
&&$\s_1$&&$\s_2$&&$h\hskip.4cm$&&$v_1$&&$v_2$&\cr
\tablerule 
&&$1$&&$1$&&$v-v'$&&$v_*$&&$ v$&\cr
\tablerule 
&&$-1$&&$-1\hskip.27cm$&&$v'-v$&&$ v$&&$v_* $&\cr
\tablerule 
&&$1$&&$-1\hskip.27cm$&&$ v'-v$\hskip-.2cm&&$ v'_*$&&$ v'$&\cr
\tablerule 
&&$-1$&&$1$&&$ v-v'\hskip-.2cm$&&$ v'$&&$ v'_*$&\cr
\tablerule 
}} 
}
\bigskip
\noi
This results in the final expression:
$$
\eqalign{
{\Cal T_{2}^{(2,1)}}(x,v)=&
{1\over 4\pi^2}\int dv_*\int dv'\int dv_*' \;
\delta(v_*+v-v_*'-v')
\cr&
\delta\left({1\over2}(v_*\/^2+v^2-v_*'\/^2-v'\/^2)\right) \;
\widehat\phi(v'-v) \, \widehat\phi(v'-v_*) \, (f'f'_*-ff_*).
}
\Eq(4.20)
$$

\bigskip

As a conclusion for the $\Cal T_2$ terms, we have eventually 
established the (desired) equality
$$
\Cal T_{2}^{(1,2)}+ 
\Cal T_{2}^{(2,1)}=\int dv_*\int dv'\int
dv_*W(v,v_*|v',v'_*) (f'f'_*-ff_*).
\Eq(4.21)
$$
This ends up the analysis of the two-particle terms.

\bigskip\bigskip


\heading 5. Three-particle terms: permutations with a fixed element 
\endheading

\numsec= 5
\reset
\numfor= 1
\bigskip


In this section we analyze $\Cal W^\pi_{3,\e}$
and $\Cal V^\pi_{3,\e}$ for the permutations 
$\pi$ with a fixed element. To simplify the notation we set 
$$
\Cal W^0_{3,\e}\quad \text{and }
\Cal V^0_{3,\e}\quad\text{ for } \pi=(1,2,3),
$$
and
$$
\Cal W^i_{3,\e}\quad \text{and }
\Cal V^i_{3,\e}, \quad i=1,2,3,
$$
for the three permutations leaving $i$ fixed.
To state the result briefly, let us readily say that the factors
$$
\Cal W^0_{3,\e},\quad
\Cal V^0_{3,\e},\quad
\text{ together with }
\Cal W^2_{3,\e},\quad
\Cal V^1_{3,\e},
$$
give a vanishing contribution.
Also, the sum
$$
\Cal W^3_{3,\e}+\Cal V^3_{3,\e}
$$
is shown to vanish asymptotically, while each of these two terms is $O(\e^{-1})$
separately. Here, anti-symmetry will play a central role.
Finally, the two terms
$$
\Cal W^1_{3,\e},\quad
\Cal V^2_{3,\e},
$$
do contribute to the limiting U-U equation through the cubic 
term. They build up the transition kernel
$\widehat \phi(v'-v)^2+ \widehat \phi(v'-v_*)^2$.
The missing cross term 
$2 \widehat \phi(v'-v) \, \widehat \phi(v'-v_*)$ in $W(v,v_*|v',v'_*)$
will come up in the
next section.

\bigskip

Let us show first that $\Cal W^0_{3,\e}$ and $\Cal V^0_{3,\e}$ 
are vanishing.
From \equ(3.19), scaling $h_1$ and $h_2$ and summing on $\s_1,\s_2$,
we have
$$
\eqalign{
\widehat{\Cal W}^0_{3,\e}&=
 +  4{\e^{-1}\over (2\pi)^6}
\int_0^{t_1} dt_2\int dh_1\int dh_2 \;
\widehat\phi(\e h_1) \, \widehat\phi(\e h_2)
\cr& 
\quad
\sin\left(\e 
{h_1\cdot \eta\over 2}
\right)
\sin\left(\e {h_2\over 2} \cdot (\eta+(t_1-t_2)(\x-h_1)\right)
\cr&
\quad
\widehat f_0(\x-(h_1+h_2), \eta+t_1\x-(t_1h_1+t_2h_2)) \,
\widehat f_0(h_1,t_1 h_1) \,
\widehat f_0(h_2,t_2 h_2)
\cr&
=O(\e),
}
\Eq(5.1)
$$
due to the decay properties of $\widehat f_0$. The same argument 
easily leads to
$
\widehat{\Cal V}^0_{3,\e}
=O(\e).
$

\bigskip

We now pass to the computation of $\widehat{\Cal W}^1_{3,\e}$. 
This term is associated with the permutation 
$\pi=(1,3,2)$. Upon Fourier
transforming in $x$ the relation \equ(til) for $\widetilde f^\pi_j$ 
(with $j=3$),
and using the change of
variables $y_1=(x_2+x_3)/2$, $y_2=(x_2-x_3)/\e$ in the 
corresponding formula,
we recover 
$$
\eqalign{
\widehat f^{(1,3,2)}_{3,\e}(\xi_1,&\xi_2,\xi_3,\h_1,\h_2,\h_3)=
\e^3\widehat f_0(\x_1,\eta_1)
\int dy_1\int dy_2 \;
\ex^{-i\x_2\cdot(y_1+{\e\over 2}y_2)} \,
\ex^{-i\x_3\cdot(y_1-{\e\over 2}y_2)}
\cr&
\widetilde f_0\left(y_1-\e{\eta_2-\h_3\over 4},
-y_2+{\h_2+\h_3\over  2}\right) \,
\widetilde f_0\left(y_1+\e{\eta_2-\h_3\over 4},
y_2+{\h_2+\h_3\over 2}\right).}
\Eq(5.2)
$$
Then, inserting \equ(5.2) in the formula \equ(3.19) relating the value of 
$\Cal W^\pi_{3,\e}$, we arrive at  
$$
\eqalign{
\widehat{\Cal W}_{3,\e}^1(\x,\h)
&=
-{\e^{-4}\over (2\pi)^6}
\sum_{\s_1,\s_2=\pm 1}\s_1\s_2 \int_0^{t_1} 
dt_2 \int dh_1 \int dh_2 \int dy_1\int  dy_2 \;
\widehat\phi(h_1) \, \widehat\phi(h_2)
\cr&
\quad
\ex^{-i\left({h_1+h_2\over \e}\cdot y_1+{h_1-h_2\over 2}\cdot y_2\right)}
\,
\ex^{-i{\s_1\over 2}h_1\cdot \eta} \,
\ex^{-i{\s_2\over 2}h_2\cdot\left(\h+(t_1-t_2)(\x-{h_1\over \e}) \right)}
\cr&
\quad
\widehat f_0
\left(\x-{h_1+h_2\over \e},\h+t_1\x-{t_1h_1+t_2h_2\over \e}\right)
\cr&
\quad
\widetilde f_0
\left(y_1-{t_1h_1-t_2h_2\over 4}, {t_1h_1+t_2h_2\over 2\e}-y_2\right)
\cr&
\quad
\widetilde f_0
\left(y_1+{t_1h_1-t_2h_2\over 4}, {t_1h_1+t_2h_2\over 2\e}+y_2\right).
}
\Eq(5.3)
$$
Changing variables
$$
k={h_1+h_2\over \e}, \quad h=h_2, \quad s={t_1-t_2\over \e},
$$
we obtain
$$
\eqalign{&
\widehat{\Cal W}_{3,\e}^1(\x,\h)=
-{1\over (2\pi)^6}
\sum_{\s_1,\s_2=\pm 1}\s_1\s_2
\int_0^{t_1\over \e} ds \int dh \int dk \int dy_1\int dy_2 \;
\widehat\phi(-h+\e k) \, \widehat\phi(h)
\cr&
\ex^{-i{\s_1\over 2} \eta \cdot (-h + \e k)} \,
\ex^{-i{\s_2\over 2}h\cdot(\h+sh +\e s(\x-k))} \,
\ex^{-ik\cdot y_1-i y_2\left(-h+\e {k \over 2}\right)} \,
\widehat f_0(\x-k,\h+t_1(\x-k) +sh)
\cr&
\widetilde f_0
\left(y_1+{t_1h\over 2}-\e{kt_1+sh\over 4}, {kt_1-sh\over 2}-y_2\right)
\,
\widetilde f_0
\left(y_1-{t_1h\over 2}+\e{kt_1+sh\over 4}, {kt_1-sh\over 2}+y_2\right).
}
\Eq(5.4)
$$
We are now in position to identify 
the rigorous limit of $\widehat{\Cal W}^1_{3,\e}$,
using the assumed decay of $\widehat\phi$ and $\widehat f_0$. 
The argument are those used
in the previous section: we do not repeat them here and in the sequel.
Passing to the limit we get, eventually,
$$
\eqalign{
\widehat{\Cal W}&_{3}^1(\x,\h)=
-{1\over (2\pi)^6}
\sum_{\s_1,\s_2=\pm 1}\s_1\s_2
\int_0^{+\infty} ds \int dh \int dk \int dy_1\int  dy_2 \;
|\widehat\phi(h)|^2
\cr&
\ex^{i{\s_1-\s_2\over 2}h\cdot \eta}\ex^{-i{\s_2\over 2}h^2s}\,
\ex^{-ik\cdot y_1+ih\cdot y_2} \,
\widehat f_0(\x-k,\h+t_1(\x-k) +sh)
\cr&
\widetilde f_0\left(y_1+{t_1h\over 2}, {kt_1-sh\over 2}-y_2\right) \,
\widetilde f_0\left(y_1-{t_1h\over 2}, {kt_1-sh\over 2}+y_2\right).
}
\Eq(5.5)
$$
Last, we go back to the $(x,v)$ variables,
computing the inverse Fourier 
transform of the above term. This gives 
$$
\eqalign{
\Cal W_{3}^1(x,v)&=
-\pi
\sum_{\s_1,\s_2=\pm 1}\s_1\s_2
\int dv_1\int dv_2\int dv_3 \;
|\widehat\phi(v_2-v_3)|^2
\cr&
\d\left(
[v_2-v_3]
\cdot
\left[-v+v_2{1+\s_1\over 2}+v_3{1-\s_1\over 2}\right]
\right)
\,
\d\left(v-v_1+{\s_1-\s_2\over 2}(v_3-v_2)\right)
\cr&
f_0(x-v_1t_1,v_1) \,
f_0(x-v_2 t_1,v_2) \,
f_0(x-v_3t_1, v_3).
}
\Eq(5.6)
$$
This is our final expression of
$\Cal W_{3}^1$. It will be interpreted later
in terms of the $v$, $v_*$, $v'$, $v'_*$ variables of the U-U equation.

\bigskip

In a similar fashion we compute 
$\Cal V_{3,\e}^2$ and its limit $\Cal V_{3}^2$.
We write
$$
\eqalign{
\widehat 
f^{(3,2,1)}_{3,\e}(\xi_1,&\xi_2,\xi_3,\h_1,\h_2,\h_3)=
\e^3\widehat  f_0(\x_2,\eta_2)
\int dy_1\int dy_2 \;
\ex^{-i\x_1\cdot(y_1+{\e\over 2}y_2)} \,
\ex^{-i\x_3\cdot(y_1-{\e\over 2}y_2)}
\cr&
\widetilde f_0
\left(y_1-\e{\eta_1-\h_3\over 4},
-y_2+{\h_2+\h_1\over 2}\right)\,
\widetilde f_0
\left(y_1+\e{\eta_1-\h_3\over 4},
y_2+{\h_1+\h_3\over 2}\right).
}
\Eq(5.2')
$$
We insert this expression into \equ(3.20), and 
perform the change of
variables $h=h_2$, $k=(h_1-h_2)/\e$,
and $s=(t_1-t_2)/\e$. Passing to the limit 
$\e\to 0$ at once, gives the asymptotic value
$$
\eqalign{&
\widehat{\Cal V}_{3}^2(\x,\h)
=-{1\over (2\pi)^6}\sum_{\s_1,\s_2=\pm 1}\s_1\s_2
\int_0^{+\infty} ds \int dh \int dk \int dy_1\int dy_2 \; 
|\widehat\phi(h)|^2
\cr&
\qquad
\widetilde f_0
\left(y_1+{t_1h\over 2}, {\h+\x t_1-kt_1-sh\over 2}-y_2\right)
\,
\widetilde f_0
\left(y_1-{t_1h\over 2}, {\h+\x t_1-kt_1-sh\over 2}+y_2\right)
\cr&
\qquad
\widehat f_0(k,t_1k+sh) \,
\ex^{-i(\x-k)\cdot y_1} \,
\ex^{ih\cdot y_2}\ex^{-i {\s_1  \over 2} h\cdot \eta} \,
\ex^{-i{\s_2\over 2}h^2s},
}
\Eq(5.7)
$$
whose inverse Fourier transform is 
$$
\eqalign{&
\Cal V_{3}^2(x,v)
=-\pi\sum_{\s_1,\s_2=\pm 1}\s_1\s_2
\int dv_1\int dv_2\int dv_3\;
|\widehat\phi(v_1-v_3)|^2
\cr&
\qquad
\d\left([v_1-v_3]\cdot
\left[v_2-v_1{1+\s_2\over 2}-v_3{1-\s_2\over 2}\right]\right)\,
\d\left(v-v_1 {1-\s_1\over 2} -v_3{1+\s_1\over 2}\right)
\cr&
\qquad
f_0(x-v_1t_1,v_1)f_0(x-v_2 t_1,v_2) \,
f_0(x-v_3t_1, v_3).
}
\Eq(5.8)
$$

\bigskip

Before coming to the computation
of the other $\Cal W^i_3$'s and
$\Cal V^i_3$'s, let us now identify the link between the obtained
values of $\Cal W^1_3$, $\Cal V^2_3$, and the U-U equation.
The following changes of variables in $\Cal W^1_3$
\bigskip
\centerline{
\vbox{\tabskip=0pt \offinterlineskip
\def\tablerule{\noalign{\hrule}}
\halign to250pt{
\strut#& \vrule#\tabskip=1em plus2em&
\hfil#&\vrule#& \hfil#\hfil& \vrule#&
\hfil#&\vrule#& \hfil#\hfil& \vrule#&
\hfil#& \vrule#
\tabskip=0pt
\cr
\tablerule
&&$\s_1$&&$\s_2$&&$v_1$&&$v_2$&&$v_3$&\cr
\tablerule 
&&$1$&&$1$&&$v\hskip.12cm$&&$ v'_*$&&$ v_*$&\cr
\tablerule 
&&$-1$&&$\hskip-.27cm-1$&&$v\hskip.12cm$&&$ v_*$&&$ v_*'$&\cr
\tablerule 
&&$1$&&$\hskip-.27cm-1$&&$v'$&&$ v'_*$&&$ v_*$&\cr
\tablerule 
&&$-1$&&$1$&&$v'$&&$ v_*$&&$ v_*'$&\cr
\tablerule 
}} 
}
\bigskip
\noi
yields the explicit value
$$
\eqalign{
\Cal W^1_3(x,v)&=
 2  \pi \int dv_*\int dv'\int dv_*'\;
|\widehat\phi(v'-v)|^2 
\,
\cr&
\d(v+v_*-v'-v_*')
\,
\d\left({1\over 2}(v^2+v_*^2-v'\/^2-v_*'\/^2)\right) \,
[f_*f'f_*'-ff_*f'_*].}
\Eq(5.9)
$$
For $\Cal V_3^2$ we set
\bigskip
\centerline{
\vbox{\tabskip=0pt \offinterlineskip
\def\tablerule{\noalign{\hrule}}
\halign to250pt{
\strut#& \vrule#\tabskip=1em plus2em&
\hfil#&\vrule#& \hfil#\hfil& \vrule#&
\hfil#&\vrule#& \hfil#\hfil& \vrule#&
\hfil#& \vrule#
\tabskip=0pt
\cr
\tablerule
&&$\s_1$&&$\s_2$&&$v_1$&&$v_2$&&$v_3$&\cr
\tablerule 
&&$1$&&$1$&&$ v'_*\hskip.12cm$&&$ v_*$&&$ v$&\cr
\tablerule 
&&$-1$&&$\hskip-.27cm-1$&&$ v$&&$v_*$&&$ v'_*\hskip.12cm$&\cr
\tablerule 
&&$1$&&$\hskip-.27cm-1$&&$ v'_*\hskip.12cm$&&$v'$&&$ v$&\cr
\tablerule 
&&$-1$&&$1$&&$ v$&&$v'$&&$ v'_*\hskip.12cm$&\cr
\tablerule 
}} 
}
\bigskip
\noi
and obtain the final
$$
\eqalign{
\Cal V^2_3(x,v)&=
 2  \pi \int dv_*\int dv'\int dv_*' \; 
|\widehat\phi(v-v'_*)|^2 \,
\cr&
\d(v+v_*-v'-v_*') \,
\d\left({1\over 2}(v^2+v_*^2-v'\/^2-v_*'\/^2)\right)\,
[ff'f_*'-f_*'ff_*].
}
\Eq(5.10)
$$
Last, 
using the symmetry $v'\leftrightarrow v_*'$ (note that
$v'-v=-(v_*'-v_*)$ and
$v'-v_*=v-v_*'$) we finally conclude after some computations that:
$$
\eqalign{&
(\Cal W^1_3+\Cal V_3^2)(x,v)=
\pi \int dv_*\int dv'\int dv_*' \;
[(f+f_*)f'f_*'-(f'+f_*')ff_*)]
\cr&
\qquad
\d(v+v_*-v'-v_*') \,
\d\left({1\over2}(v^2+v_*^2-v'\/^2-v_*'\/^2)\right) \,
(|\widehat\phi(v'-v)|^2+|\widehat\phi(v'-v_*)|^2).}
\Eq(5.11)
$$
This 
is the desired cubic term in the U-U equation, up to the fact that we
only recover here part of the total cross-section
$W(v,v_*|v',v'_*)=[\widehat \phi(v-v')+\widehat \phi(v-v_*)]^2$. The
missing cross term will come up in the next section.

\bigskip

We now show that all other terms associated to permutations with a 
fixed element, namely
$\Cal W_{3,\e}^2$,
$\Cal V_{3,\e}^1$, and 
$\Cal W_{3,\e}^3+\Cal V_{3,\e}^3$,
give a vanishing contribution in the limit $\e\to 0$. 

\bigskip

We begin with $\Cal W_{3,\e}^2$. Inserting \equ(5.2') into \equ(3.19)
and changing variables $k=h_1\e^{-1}$, $h_2=h$, we readily obtain:
$$
\eqalign{
\widehat{\Cal W}^2_{3,\e}(\x,\h)=&
-{\e^{-1}\over (2\pi)^6}
\sum_{\s_1,\s_2=\pm1}\s_1\s_2
\int_0^{t_1}dt_2\int dk\int dh\int dy_1\int dy_2 \;
\widehat\phi(h) \, \widehat\phi(\e k)
\cr&
\ex^{-{i\over 2}(\s_1\h\cdot k\e+\s_2h\cdot(\h+(t_1-t_2)(\x-k))}
\,
\ex^{iy_2\cdot h}\,
\ex^{-i\left(y_1+\e {y_2\over 2}\right)\cdot(\x-k)} \,
\cr&
\widetilde f_0
\left(y_1-{\e\over 4}\left[\h+\x t_1-t_1k-{2\over \e}t_2 h\right];
{1\over2}(\h+\x t_1-t_1k)-y_2\right)
\cr&
\widetilde f_0
\left(y_1+{\e\over 4}\left[\h+\x t_1-t_1k-{2\over \e}t_2 h\right];
{1\over2}(\h+\x t_1-t_1k)+y_2\right)
\cr&
\widehat f_0(k, t_1 k).
}
\Eq(5.12)
$$
Summing on $\s_1$ and $\s_2$ allows to compute the limit which is:
$$
\eqalign{
\widehat{\Cal W}^2_{3}&(\x,\h)=
{4\over (2\pi)^6}\int_0^{t_1}dt_2\int dk\int dh
\int dy_1\int dy_2 \;
\widehat\phi(h)\, \widehat\phi(0)\,
\ex^{iy_2\cdot h} \,\ex^{-iy_1\cdot(\x-k)}
\cr&
{\h\cdot k \over 2}\,
\sin\left({1\over 2}h\cdot(\h+(t_1-t_2)(\x-k)\right)
\,
\widetilde f_0
\left(y_1+{1\over 2}t_2 h;{1\over 2}(\h+\x t_1-t_1k)-y_2\right)
\cr&
\widetilde f_0
\left(y_1-{1\over 2}t_2 h;{1\over2}(\h+\x t_1-t_1k)+y_2\right)
\,
\widehat f_0(k, t_1 k).
}
\Eq(5.13)
$$
Note that the product of the two $\widetilde f_0$'s is invariant under 
the change of variables $h\to -h$, $y_2\to -y_2$, as well
as the oscillatory factor $\ex^{i y_2 \cdot h}$. Therefore, 
using the 
parity of $\widehat\phi$, it follows that
$\widehat{\Cal W}^2_{3}(\x,\h)=-\widehat{\Cal W}^2_{3}(\x,\h)$.
Hence
$$
\widehat{\Cal W}^2_{3}=0.
$$

The term $\widehat{\Cal V}^1_{3}$ 
is treated in the analogous way.
We insert \equ(5.2) into \equ(3.20) and make the change of variables
$h=h_2$, $k=\e^{-1} h_1$, obtaining after summation over $\s_1$ and $\s_2$:
$$
\eqalign{ 
\widehat{\Cal V}^1_{3,\e}&(\x,\h)=
{4 \e^{-1}\over (2\pi)^6}
\int_0^{t_1}dt_2\int dk\int dh \int dy_1\int dy_2 \; 
\widehat\phi(h) \, \widehat\phi(\e k)
\cr&
\ex^{-iy_1\cdot k} \, 
\ex^{iy_2\cdot \left(h-\e {k \over 2}\right)} \,
\sin\left({\e\over 2}\h\cdot k\right) \,
\sin\left({1\over 2}h\cdot k(t_1-t_2)\right)
\cr&
\widetilde f_0
\left(y_1-{\e\over 4}\left[t_1k-{2\over \e}t_2 h\right];{1\over 2}t_1k-y_2\right)
\,
\widetilde f_0
\left(y_1+{\e\over 4}\left[t_1k-{2\over \e}t_2 h\right];{1\over 2}t_1k+y_2\right)
\cr&
\widehat f_0(\x-k, \h+t_1 (\x-k)).
}\Eq(5.14)
$$
This term clearly goes to
$$
\eqalign{
\widehat{\Cal V}^1_{3}&(\x,\h)=
{4\over (2\pi)^6}
\int_0^{t_1}dt_2\int dk\int dh \int dy_1\int dy_2 \;
\widehat\phi(h)\widehat\phi(0)
\cr&
\quad
\ex^{-iy_1\cdot k} \,
\ex^{iy_2\cdot h}\,
\left({1\over 2}\h\cdot k\right)\,
\sin\left({1\over 2}h\cdot k(t_1-t_2)\right)
\,
\widetilde f_0
\left(y_1+{1\over 2}t_2 h;{1\over 2}t_1k-y_2\right)
\cr&
\quad
\widetilde f_0
\left(y_1-{1\over 2}t_2 h;{1\over 2}t_1 k+y_2\right)
\,
\widehat f_0(\x-k, \h+t_1 (\x-k)).
}
\Eq(5.15)
$$
Hence $\widehat{\Cal V}^1_{3}=0$ for the same reason as before.

\bigskip

To end up this paragraph, let us last prove that the sum
$\Cal W^3_{3,\e}+\Cal V_{3,\e}^3$ vanishes asymptotically.
First, we write
$$
\eqalign{
\widehat f^3_{3,\e}(\xi_1,&\xi_2,\xi_3,\h_1,\h_2,\h_3)=
\e^3\widehat f_0(\x_3,\eta_3)
\int dy_1\int dy_2 \;
\ex^{-i\x_1\cdot(y_1+{\e\over 2}y_2)} \,
\ex^{-i\x_2\cdot(y_1-{\e\over2}y_2)}
\cr&
\widetilde f_0
\left(y_1-\e{\eta_1-\h_2\over 4}, -y_2+{\h_1+\h_2\over  2}\right)
\,
\widetilde f_0
\left(y_1+\e{\eta_1-\h_2\over 4}, y_2+{\h_1+\h_2\over 2}\right).
}
\Eq(5.2'')
$$
We insert this formmula
in \equ(3.19),
and perform the change of
variables $h_1=h$, $h_2=\e k$.
In this way we recover, upon computing the
sum $\sum \s_1 \s_2$,
$$
\eqalign{
\widehat{\Cal W}^3_{3,\e}&(\x,\h)=
{4\e^{-1}\over (2\pi)^6}
\int_0^{t_1}dt_2\int dk\int dh\int dy_1\int dy_2\;
\widehat\phi(h) \, \widehat\phi(\e k)
\cr&
\ex^{iy_2\cdot h}\,
\ex^{-i(y_1+{\e \over 2}y_2)\cdot(\x- k)}\,
\sin\left({1\over 2}\h\cdot h\right) \,
\sin\left({\e k\over 2}\cdot (\h+(t_1-t_2)\x)
-{1\over 2}k\cdot h(t_1-t_2)\right)
\cr&
\widetilde f_0
\left(y_1-{\e\over 4}\left[\h+\x t_1-kt_2-{2\over \e}t_1 h\right];
{1\over2}(h+\x t_1-t_2k)-y_2\right)
\cr&
\widetilde f_0
\left(y_1+{\e\over 4}\left[\h+\x t_1-kt_2-{2\over \e}t_1 h\right];
{1\over2}(\h+\x t_1-t_2k)+y_2\right)
\cr&
\widehat f_0(k,t_2k).
}\Eq(5.16)
$$
Similarly, using \equ(5.2''), \equ(3.20), and performing again the change of
variables  $h_1=h$, $h_2=\e k$, we obtain:
$$
\eqalign{
\widehat{\Cal V}^3_{3,\e}(\x,\h)=
&{4\e^{-1}\over (2\pi)^6}
\int_0^{t_1}dt_2\int dk\int dh\int dy_1\int dy_2\;
\widehat\phi(h)\, \widehat\phi(\e k)
\cr&
\ex^{iy_2\cdot h}\,
\ex^{-iy_1\cdot(\x-k)}\,
\ex^{-i{\e\over 2}y_2\cdot(\x+k)}
\cr&
\sin\left({1\over 2}\h\cdot h\right)\,
\sin\left({1\over 2}k\cdot h(t_1-t_2)\right)
\cr&
\widetilde f_0
\left(y_1+{\e\over 4}\left[\h+\x t_1+kt_2-{2\over \e}t_1 h\right];
{1\over2}(\h+\x t_1-t_2k)+y_2\right)
\cr&
\widetilde f_0
\left(y_1-{\e\over 4}\left[\h+\x t_1+kt_2-{2\over \e}t_1 h\right];
{1\over2}(\h+\x t_1-t_2k)-y_2\right)
\cr&
\widehat f_0(k,t_2k).
}
\Eq(5.17)
$$
Hence, both terms
$\widehat{\Cal W}^3_{3,\e}$ and $\widehat{\Cal V}^3_{3,\e}$
are $O(\e^{-1})$. However we have the following expansions:
$$
\e\widehat{\Cal W}^3_{3,\e}=A_0+\e A_1+O(\e^2),
\quad
\e\widehat{\Cal V}^3_{3,\e}=B_0+\e B_1+O(\e^2),
$$
and it is easy to realize that $A_0=-B_0$.
Moreover, after some straightforward calculation, we obtain at the next
order:
$$
\eqalign{&
A_1+B_1=
{4\e^{-1}\over (2\pi)^6}\int_0^{t_1}dt_2\int dk\int dh
\int dy_1\int dy_2\;
\widehat\phi(h)\,
\widehat\phi(0)
\ex^{iy_2\cdot h}\, \ex^{-iy_1\cdot(\x-k)}
\,
\sin\left({1\over 2}\h\cdot h\right)
\cr&
\widehat f_0(k,t_2k)
\widetilde f_0
\left(y_1-{1\over 2}t_1 h;{1\over2}(\h+\x t_1-t_2k)+y_2\right)\,
\widetilde f_0
\left(y_1+{1\over 2}t_1 h;{1\over2}(\h+\x t_1-t_2k)-y_2\right)
\cr&
\left\{
{k\over 2}\cdot(\h+(t_1-t_2)\x) \,
\cos\left({1\over 2}k\cdot h(t_1-t_2)\right)+
\right.
\cr&
\sin\left({1\over2}k\cdot h(t_1-t_2)\right)
\left.
\left({1\over 2} t_2k\cdot \n_x
\log
{\widetilde f_0((y_1-{1\over 2}t_1 h;{1\over2}(\h+\x t_1-t_2k)+y_2)
\over
\widetilde f_0((y_1+{1\over 2}t_1 h;{1\over 2}(\h+\x t_1-t_2k)-y_2)}
-iy_2\cdot k\right)\right\}.
}
\Eq(5.18)
$$
Let us now exchange $h\to -h$ and $y_2\to -y_2$. The term in braces is 
invariant  because the
$\log$ term changes its sign. All the other terms are invariant but
$\sin(\h\cdot h/2)$, which changes sign. Therefore
$A_1+B_1=-(A_1+B_1),$
hence $A_1+B_1=0$. This shows that $\Cal W_{3,\e}^3+\Cal V_{3,\e}^3$ 
vanishes in the
limit $\e\to 0$.

\bigskip\bigskip


\heading 6. Three-particle terms: cyclic permutations\endheading

\numsec=6
\reset
\numfor= 1
\bigskip


We still have to evaluate $\Cal W^\pi_{3,\e}$, $\Cal
W^{\pi^{-1}}_{3,\e}$,
$\Cal V^\pi_{3,\e}$ and $\Cal V^{\pi^{-1}}_{3,\e}$ with $\pi=(2,3,1)$
and
$\pi^{-1}=(3,1,2)$.

\bigskip

We first observe, for later convenience, the two relations
$$
\eqalign{
\widehat f_3^{\pi^{-1}}(\x_1,\x_2,\x_3;&\h_1,\h_2,\h_3)=
\int dx_1 \int dx_2 \int dx_3 \;
\ex^{-i\sum_{k=1}^3 x_k\cdot\x_k}
\cr&
\widetilde f_0
\left({x_1+x_2\over 2}+{\e\over 4}(\h_1-\h_2);
{x_1-x_2\over\e}+{\h_1+\h_2\over 2}\right)
\cr&
\widetilde f_0
\left({x_2+x_3\over 2}+{\e\over 4}(\h_2-\h_3);
{x_2-x_3\over \e}+{\h_2+\h_3\over 2}\right)
\cr&
\widetilde f_0
\left({x_3+x_1\over 2}+{\e\over 4}(\h_3-\h_1);
{x_3-x_1\over  \e}+{\h_3+\h_1\over 2}\right),
}
\Eq(6.1)
$$
\vskip -0.2cm
$$
\eqalign{
\widehat  f_3^{\pi}(\x_1,\x_2,\x_3;&\h_1,\h_2,\h_3)=
\int dx_1\int dx_2\int dx_3 \;
\ex^{-i\sum_{k=1}^3 x_k\cdot\x_k}
\cr&
\widetilde f_0
\left({x_2+x_1\over 2}+{\e\over 4}(\h_2-\h_1);
{x_2-x_1\over  \e}+{\h_2+\h_1\over 2}\right)
\cr&
\widetilde f_0
\left({x_3+x_2\over 2}+{\e\over 4}(\h_3-\h_2);
{x_3-x_2\over  \e}+{\h_3+\h_2\over 2}\right)
\cr&
\widetilde f_0
\left({x_1+x_3\over 2}+{\e\over 4}(\h_1-\h_3);
{x_1-x_3\over \e}+{\h_1+\h_3\over 2}\right).
}
\Eq(6.2)
$$

\bigskip

Armed with these expressions, we begin with the computation of 
$\widehat{\Cal W}_{3,\e}^{\pi^{-1}}$ and
$\widehat{\Cal W}_{3,\e}^{\pi}$.
To do so, we insert \equ(6.1) and \equ(6.2) in 
the general formula 
\equ(3.19) relating the value of the $\widehat{\Cal W}_{3,\e}^{\pi}$'s.
In the so-obtained formulae,
we also change variables,
$h_1\to  -h_1$, $h_2\to -h_2$ for $\pi^{-1}$, and
$\s_1\to -\s_1$, $\s_2\to -\s_2$ for $\pi$. 
With this new set of variables, the $\xi$'s and $\eta$'s
involved in \equ(3.19) are
$$
\eqalign{
&\x_1=\x\pm{h_1+h_2\over \e},\quad\hskip 1.5cm
\x_2=\mp {h_1\over \e},\quad\hskip.3cm
\x_3=\mp {h_2\over\e},
\cr&
\h_1=\h+\x t_1\pm {t_1h_1+t_2h_2\over \e},\quad
\h_2=\mp {t_1h_1\over \e},\quad
\h_3=\mp {t_2h_2\over \e},}
\Eq(6.3)
$$
for $\pi^{-1}$ and $\pi$ respectively.
Also,
the phases appearing in \equ(3.19) are:
$$
\eqalign{
&
\Cal S_{\pi^-1}=
{\s_1\over 2}h_1\cdot \h
+
{\s_2\over 2}h_2\left(\h+(t_1-t_2)\left[\x+{h_1\over\e}\right]\right),
\cr&
\Cal S_{\pi}=
{\s_1\over 2}h_1\cdot \h
+
{\s_2\over2}h_2\left(\h+(t_1-t_2)\left[\x-{h_1\over\e}\right]\right).
}
\Eq(6.5)
$$
All this 
gives in \equ(3.19), the two values
$$
\eqalign{
\widehat{\Cal W}_{3,\e}^{\pi^{-1}}(\x,\h)=&
-{\e^{-7}\over (2\pi)^6}
\sum_{\s_1\s_2=\pm 1}\s_1\s_2
\int_0^{t_1}dt_2\int dh_1\int dh_2\int dx_1\int dx_2\int dx_3\;
\cr&
\widetilde f_0
\left({x_1+x_2\over 2}+{2 t_1h_1+t_2h_2\over 4}+{\e\over 4}\bar\h;
{x_1-x_2\over\e}+{t_2h_2\over 2\e}+{\bar\h\over 2}
\right)
\cr&
\widetilde f_0
\left({x_2+x_3\over 2}-{t_1h_1-t_2h_2\over 4};
{x_2-x_3\over \e}-{t_1h_1+t_2h_2\over 2\e}
\right)
\cr&
\widetilde f_0
\left({x_3+x_1\over 2}-{t_1h_1+2t_2h_2\over 4}-{\e\over 4}\bar\h;
{x_3-x_1\over\e}+{t_1h_1\over 2\e}+{\bar\h\over 2}
\right)
\cr&
\widehat\phi(h_1) \, \widehat\phi(h_2) \,
\ex^{-ix_1\cdot (\x+{h_1+h_2\over \e})}\,
\ex^{i x_2\cdot{h_1\over\e}}\,
\ex^{ix_3\cdot {h_2\over \e}}\,
\ex^{i\Cal S_{\pi^{-1}}},
}
\Eq(6.6)
$$
\vskip -0.25cm
$$
\eqalign{
\widehat{\Cal W}_{3,\e}^{\pi}(\x,\h)=&
-{\e^{-7}\over (2\pi)^6}
\sum_{\s_1\s_2=\pm1}\s_1\s_2
\int_0^{t_1}dt_2\int dh_1\int dh_2\int dx_1\int dx_2\int dx_3\;
\cr&
\widetilde f_0
\left({x_2+x_1\over 2}+{2 t_1h_1+t_2h_2\over 4}-{\e\over 4}\bar\h;
{x_2-x_1\over\e}-{t_2h_2\over 2\e}+{\bar\h\over 2}
\right)
\cr&
\widetilde f_0
\left({x_3+x_2\over 2}-{t_1h_1-t_2h_2\over 4};
{x_3-x_2\over\e}+{t_1h_1+t_2h_2\over 2\e}
\right)
\cr&
\widetilde f_0
\left({x_1+x_3\over 2}-{t_1h_1+2t_2h_2\over 4}+{\e\over 4}\bar\h;
{x_1-x_3\over\e}-{t_1h_1\over 2\e}+{\bar\h\over 2}
\right)
\cr&
\widehat\phi(h_1)\, \widehat\phi(h_2)\,
\ex^{-ix_1\cdot (\x-{h_1+h_2\over \e})}\,
\ex^{-i x_2\cdot{h_1\over\e}}\,
\ex^{-ix_3\cdot {h_2\over \e}}\,
\ex^{i\Cal  S_{\pi}},
}
\Eq(6.7)
$$
where we use the notation
$
\bar\h=\h+\x t_1.
$

Now, we perform the following natural change of variables: 
$$
\eqalign{
&
\text{for } \pi^{-1}:
\qquad
x_2= x_1+{1\over 2}t_2 h_2-\e y_1,
\quad
x_3=
x_1-{1\over2}t_1h_1+\e y_3,
\cr&
\text{for } \pi:
\qquad
\quad
x_2= x_1+{1\over 2}t_2 h_2+\e y_1,
\quad
x_3=x_1-{1\over 2}t_1h_1-\e y_3.
}
\Eq(6.9)
$$
In both cases 
$x_1$ is unchanged.
This finally gives the two relations
$$
\eqalign{
\widehat{\Cal W}_{3,\e}^{\pi^{-1}}(\x,\h)=&
-{\e^{-1}\over  (2\pi)^6}
\sum_{\s_1\s_2=\pm 1}\s_1\s_2
\int_0^{t_1}dt_2\int dh_1\int dh_2\int dx_1\int  dy_1 
\int dy_3
\cr&
\widetilde f_0\left(
x_1+{t_1h_1+t_2h_2\over 2}
+{\e\over 2}\left[{\bar\h\over 2}-y_1\right]; 
y_1+{\bar\h\over 2}
\right)\cr&
\widetilde f_0
\left(x_1-{t_1h_1-t_2h_2\over 2}
-{\e\over 2}(y_1-y_3); -y_1-y_3\right)
\cr&
\widetilde f_0
\left(x_1-{t_1h_1+t_2h_2\over 2}
-{\e\over 2}\left[{\bar\h\over 2}-y_3\right];
y_3+{\bar\h\over 2}\right)
\cr&
\widehat\phi(h_1)\,
\widehat\phi(h_2)\, 
\ex^{-ix_1\cdot \x}\,
\ex^{-iy_1\cdot h_1}\,
\ex^{iy_3\cdot h_2}\,
\ex^{-{i\over 2\e}h_1\cdot h_2(t_1-t_2)}\,
\ex^{i\Cal S_{\pi^{-1}}}, 
}
\Eq(6.10)
$$
$$
\eqalign{
\widehat{\Cal W}_{3,\e}^{\pi}(\x,\h)=&
-{\e^{-1}\over (2\pi)^6}
\sum_{\s_1\s_2=\pm 1}\s_1\s_2
\int_0^{t_1}dt_2\int dh_1\int dh_2\int dx_1\int dy_1
\int dy_3
\cr&
\widetilde f_0
\left(
x_1+{t_1h_1+t_2h_2\over 2}
-{\e\over 2}\left[{\bar\h\over 2}-y_1\right];
y_1+{\bar\h\over 2}\right)
\cr&
\widetilde f_0
\left(x_1-{t_1h_1-t_2h_2\over 2} + {\e\over 2}(y_1-y_3); 
-y_1-y_3\right)
\cr&
\widetilde f_0
\left(x_1-{t_1h_1+t_2h_2\over 2}
+{\e\over 2}\left[{\bar\h\over 2}-y_3\right];
y_3+{\bar\h\over 2}\right)
\cr&
\widehat\phi(h_1) \, \widehat\phi(h_2)\,
\ex^{-ix_1\cdot \x}\,
\ex^{-iy_1\cdot h_1}\,
\ex^{iy_3\cdot h_2}\,
\ex^{{i\over 2\e}h_1\cdot h_2(t_1-t_2)}\,
\ex^{i\Cal S_{\pi}}.
}
\Eq(6.11)
$$

\bigskip

Next, we come to the computation of
$\widehat{\Cal V}^{\pi}_{3,\e}$ and $\widehat{\Cal V}^{\pi^{-1}}_{3,\e}$.
We insert
\equ(6.1) and
\equ(6.2) in \equ(3.20). We also change 
$h_1\to -h_1$, $\s_2 \to -\s_2$ 
for $\pi^{-1}$ and 
$h_2\to -h_2$, $\s_2\to-\s_2$ for $\pi$.
This gives the two identities
$$
\eqalign{
\widehat{\Cal V}_{3,\e}^{\pi^{-1}}(\x,\h)=&
-{\e^{-7}\over (2\pi)^6}
\sum_{\s_1\s_2=\pm 1}\s_1\s_2
\int_0^{t_1}dt_2\int dh_1\int dh_2\int dx_1\int dx_2\int dx_3
\cr&
\widetilde f_0
\left({x_1+x_2\over 2}+{2 t_1h_1+t_2h_2\over 4}+{\e\over 4}\bar\h;
{x_1-x_2\over \e}-{t_2h_2\over 2\e}+{\bar\h\over 2}\right)
\cr&
\widetilde f_0
\left({x_3+x_1\over 2}-{t_1h_1-t_2h_2\over 4}-{\e\over 4}\bar\h;
{x_3-x_1\over \e}+{t_1h_1+t_2h_2\over 2\e}+{\bar\h\over 2}\right)
\cr&
\widetilde f_0
\left({x_2+x_3\over 2}-{t_1h_1+2t_2h_2\over 4};
{x_2-x_3\over \e}-{t_1h_1\over 2\e}\right)
\cr&
\widehat\phi(h_1) \, \widehat\phi(h_2) \,
\ex^{-ix_1\cdot (\x+{h_1\over \e})} \,
\ex^{i x_2\cdot{h_1+h_2\over \e}} \,
\ex^{-ix_3\cdot {h_2\over\e}} \,
\ex^{i\widetilde{\Cal S}_{\pi^{-1}}},
}
\Eq(6.12)
$$
$$
\eqalign{
\widehat{\Cal V}_{3,\e}^{\pi}(\x,\h)=&
-{\e^{-7}\over (2\pi)^6}
\sum_{\s_1\s_2=\pm 1}\s_1\s_2
\int_0^{t_1}dt_2\int dh_1\int dh_2\int  dx_1\int dx_2\int dx_3
\cr&
\widetilde f_0
\left({x_2+x_1\over 2}+{2 t_1h_1+t_2h_2\over 4}-{\e\over 4}\bar\h;
{x_2-x_1\over \e}+{t_2h_2\over 2\e}+{\bar\h\over 2}\right)
\cr&
\widetilde f_0
\left({x_1+x_3\over 2}-{t_1h_1-t_2h_2\over 4}+{\e\over 4}\bar\h;
{x_1-x_3\over \e}-{t_1h_1+t_2h_2\over 2\e}+{\bar\h\over 2}\right)
\cr&
\widetilde f_0
\left({x_3+x_2\over 2}-{t_1h_1+2t_2h_2\over 4};
{x_3-x_2\over \e}+{t_1h_1\over 2\e}\right)
\cr&
\widehat\phi(h_1) \, \widehat\phi(h_2)\,
\ex^{-ix_1\cdot (\x-{h_1\over \e})} \,
\ex^{-i x_2\cdot{h_1+h_2\over \e}}\,
\ex^{ix_3\cdot {h_2\over \e}}\,
\ex^{i\widetilde{\Cal S}_{\pi}}.
}
\Eq(6.13)
$$
Here the phases are:
$$
\eqalign{
&
\widetilde{\Cal S}_{\pi^{-1}}=
{\s_1\over 2}h_1\cdot \h
+{\s_2\over 2}h_2\cdot {h_1\over \e}(t_1-t_2),
\cr&
\widetilde{\Cal S}_{\pi}=
{\s_1\over 2}h_1\cdot \h
-{\s_2\over 2}h_2\cdot {h_1\over \e}(t_1-t_2).
}
\Eq(6.15)
$$
We make the following natural change of variable
$$
\eqalign{
&
\text{for } \pi^{-1}: \qquad
x_1=x'_1+{1\over 2}t_2h_2,\quad
x_2=x'_1-\e y_1,\quad
x_3=x'_1-{1\over 2}t_1 h_1-\e\left(y_1+y_3+{\bar\h\over2}\right),
\cr&
\text{for } \pi: \qquad
x_1=x'_1+{1\over 2}t_2h_2,\quad
x_2=x'_1+\e y_1,\quad
x_3=x'_1 - {1\over 2}t_1 h_1
+ \e\left(y_1+y_3+{\bar\h\over2}\right).
}
\Eq(6.17)
$$
With this change of variables, we eventually obtain
$$
\eqalign{
\widehat{\Cal V}_{3,\e}^{\pi^{-1}}&(\x,\h)=
-{\e^{-1}\over (2\pi)^6}
\sum_{\s_1\s_2=\pm 1}\s_1\s_2
\int_0^{t_1}dt_2\int dh_1\int dh_2\int dx'_1 \int  dy_1  \int dy_3
\cr&
\widetilde f_0
\left(
x'_1+{t_1h_1+t_2h_2\over 2}+{\e\over 2}\left[{\bar\h\over 2}-y_1\right];
y_1+{\bar\h\over 2}\right)
\cr&
\widetilde f_0
\left(
x'_1-{t_1h_1-t_2h_2\over 2}-{\e\over 2}(\bar\h+ y_1+y_3);
-y_1-y_3\right)
\cr&
\widetilde f_0
\left(
x'_1-{t_1h_1+t_2h_2\over 2}-{\e\over 2}\left[{\bar\h\over2}+2y_1+y_3\right]; 
y_3+{\bar\h\over 2}\right)
\cr&
\widehat\phi(h_1)\, \widehat\phi(h_2)\,
\ex^{-ix'_1\cdot \x}\,
\ex^{-iy_1\cdot h_1}\,
\ex^{iy_3\cdot h_2}\,
\ex^{{i\over 2\e}h_1\cdot h_2(t_1-t_2)}\,
\ex^{i{h_2\over 2}\cdot(\h +\x(t_1-t_2))}\,
\ex^{i\widetilde{\Cal S}_{\pi^{-1}}},
}
\Eq(6.18)
$$
as well as
$$
\eqalign{
\widehat{\Cal V}_{3,\e}^{\pi}(\x,\h)=&
-{\e^{-1}\over(2\pi)^6}
\sum_{\s_1\s_2=\pm 1}\s_1\s_2
\int_0^{t_1}dt_2\int dh_1\int dh_2\int dx'_1 \int dy_1  \int dy_3
\cr&
\widetilde f_0
\left(
x'_1+{t_1h_1+t_2h_2\over 2}-{\e\over 2}\left[{\bar\h\over 2}-y_1\right]; 
y_1+{\bar\h\over 2}\right)
\cr&
\widetilde f_0
\left(
x'_1-{t_1h_1-t_2h_2\over 2}+{\e\over 2}(\bar\h+ y_1+y_3);
-y_1-y_3\right)
\cr&
\widetilde f_0
\left(
x'_1-{t_1h_1+t_2h_2\over 2}+{\e\over 2}\left[{\bar\h\over2}+2y_1+y_3\right]; 
y_3+{\bar\h\over 2}\right)
\cr&
\widehat\phi(h_1)\, \widehat\phi(h_2)\,
\ex^{-ix'_1\cdot \x}\,
\ex^{-iy_1\cdot h_1}\,
\ex^{iy_3\cdot h_2}\,
\ex^{-{i\over 2\e}h_1\cdot h_2(t_1-t_2)}\,
\ex^{i{h_2\over 2}\cdot(\h+\x(t_1-t_2))}\,
\ex^{i\widetilde{\Cal S}_{\pi}}.
}
\Eq(6.19)
$$

\bigskip

Let us come to the computation of the limit of the above
four terms
$\widehat{\Cal W}^{\pi^{-1}}_{3,\e}$,
$\widehat{\Cal W}^{\pi}_{3,\e}$,
$\widehat{\Cal V}^{\pi^{-1}}_{3,\e}$,
$\widehat{\Cal V}^{\pi}_{3,\e}$.
The phases carried by these terms
are respectively:
$$
{\s_1h_1+\s_2h_2\over2}\cdot \h
+{\s_2\over 2} h_2\cdot\x(t_1-t_2)
-x_1\cdot\x-y_1\cdot h_1+y_3\cdot h_2
-{1-\s_2\over 2}{t_1-t_2\over \e}h_1\cdot h_2,
$$
$$
{\s_1h_1+\s_2h_2\over2}\cdot \h
+{\s_2\over 2} h_2\cdot\x(t_1-t_2)
-x_1\cdot\x-y_1\cdot h_1+y_3\cdot h_2
+{1-\s_2\over 2}{t_1-t_2\over \e}h_1\cdot h_2,
$$
$$
{\s_1h_1+h_2\over2}\cdot \h
+{1\over 2} h_2\cdot\x(t_1-t_2)
-x'_1\cdot\x-y_1\cdot h_1+y_3\cdot h_2
+{1+\s_2\over 2} {t_1-t_2\over \e}h_1\cdot h_2,
$$
$$
{\s_1h_1+h_2\over2}\cdot \h
+{1\over 2} h_2\cdot\x(t_1-t_2)
-x'_1\cdot\x-y_1\cdot h_1+y_3\cdot h_2
-{1+\s_2\over 2}  {t_1-t_2\over \e}h_1\cdot h_2.
$$
Denoting by 
$\widehat{\Cal W}^{\pi^{-1},\pm}_{3,\e}$,
$\widehat{\Cal W}^{\pi,\pm}_{3,\e}$,
$\widehat{\Cal V}^{\pi^{-1},\pm}_{3,\e}$
$\widehat{\Cal V}^{\pi,\pm}_{3,\e}$,
the eight terms relative to the values of
$\s_2=\pm 1$, we realize that
$\widehat{\Cal W}^{\pi^{-1},+}_{3,\e}$,
$\widehat{\Cal W}^{\pi,+}_{3,\e}$,
$\widehat{\Cal V}^{\pi^{-1},-}_{3,\e}$
$\widehat{\Cal V}^{\pi,-}_{3,\e}$,
have only slowly varying phases,
so that they are individually $O(\e^{-1})$.
However, setting
$$
\eqalign{
&
\e\widehat{\Cal W}^{\pi^{-1},+}_{3,\e}=A_0+A_1\e+O(\e^2),
\quad
\e\widehat{\Cal W}^{\pi,+}_{3,\e}=B_0+B_1\e+O(\e^2),
\cr&
\e\widehat{\Cal V}^{\pi^{-1},-}_{3,\e}=C_0+C_1\e+O(\e^2),
\qquad
\e\widehat{\Cal V}^{\pi^,-}_{3,\e}=D_0+D_1\e+O(\e^2),
}
\Eq(6.24)
$$
an easy first order Taylor expansion gives 
$A_0+C_0=B_0+D_0=0$ and $A_1+B_1=C_1+D_1=0$. 
Hence
$$
\lim_{\e\to 0}\Big(
\widehat{\Cal W}^{\pi^{-1},+}_{3,\e}
+\widehat{\Cal W}^{\pi,+}_{3,\e}
+\widehat{\Cal V}^{\pi^{-1},-}_{3,\e}
+\widehat{\Cal V}^{\pi,-}_{3,\e}
\Big)=0.
\Eq(6.25)
$$
For the other terms, which carry a rapidly oscillating phases,
it is natural to
rescale time, setting
$s=\e^{-1}(t_1-t_2)$. 
Then, an easy computation shows
$$
\lim_{\e\to 0} \Big(
\widehat{\Cal W}^{\pi,-}_{3,\e}
+
\widehat{\Cal W}^{\pi^{-1},-}_{3,\e}\Big)
=:\widehat{\Cal W}^\pi_3,
\Eq(6.26)
$$
where
$$
\eqalign{
\widehat{\Cal W}_3^\pi(\x,\h)=&
{1\over (2\pi)^6}
\sum_{\s_1\=\pm 1}\s_1
\int_0^{\infty}ds\int dh_1\int dh_2 \;
\widehat\phi(h_1)\, \widehat\phi(h_2)\,
\Big(
\ex^{-ih_1\cdot h_2 s}+\ex^{ih_1\cdot h_2 s}
\Big)
\cr&
\int dx_1\int dy_1\int dy_3\;
\ex^{{i\over 2}(\s_1h_1-h_2)\cdot \h}\,
\ex^{-ix_1\cdot\x}\,
\ex^{-y_1\cdot h_1}\, 
\ex^{i y_3\cdot h_2}
\cr&
\widetilde f_0
\left(x_1+t_1{h_1+h_2\over2}; y_1+{\bar\h\over 2}\right)
\widetilde f_0\left(x_1-t_1{h_1-h_2\over 2}; -y_1-y_3\right)
\cr&
\widetilde f_0
\left( x_1-t_1{h_1+h_2\over 2};y_3+{\bar\h\over 2}\right).
}
\Eq(6.27)
$$
Finally, taking the inverse Fourier transform of this term, we obtain:
$$
\eqalign{
\Cal W_3^\pi(x,v)=&
2\pi\int dv_1\int dv_2\int dv_3\; 
\widehat\phi(v_2-v_1)\,\widehat\phi(v_3-v_2)
\cr&
[\d(v+v_2-v_1-v_3)-\d(v-v_3)]\,
\d((v_2-v_1)\cdot(v_3-v_2))
\cr&
f_0(x-v_1t_1,v_1)f_0(x-v_2t_1,v_2)f_0(x-v_3t_1,v_3)
\cr
=&
2\pi\int dv_*\int dv_*'\int dv'\;
\widehat\phi(v'-v_*)\, 
\widehat\phi(v'-v)\, 
\big\{f_*f'f_*'-ff_*f'\big\}\cr&
\d(v+v_*-v'-v'_*)\,
\d\left({1\over2}(v^2+v_*^2-v'\/^2-v_*'\/^2)\right).}
\Eq(6.28)
$$
In the similar way we compute the limit
$$
\lim_{\e\to 0}
\left(
\Cal V_{3,\e}^{\pi^{-1},+}
+
\Cal V_{3,\e}^{\pi,+}
\right)
=
\Cal V_{3}^{\pi},
$$
whose inverse Fourier transform admits the value
$$
\eqalign{\Cal V_3^\pi(x,v)=&
2\pi\int dv_1\int dv_2\int dv_3\;
\widehat\phi(v_2-v_1)\, \widehat\phi(v_3-v_2)
\cr&
[\d( v-v_2 )-\d(v- v_1 )] \, 
\d((v_2-v_1) \cdot(v_3-v_2))
\cr&
f_0(x-v_1t_1,v_1)f_0(x-v_2t_1,v_2)f_0(x-v_3t_1,v_3)
\cr
=&
2\pi\int dv_*\int dv_*'\int dv'\;
\widehat\phi(v'-v_*)\,\widehat\phi(v'-v)
\big\{ff'f_*'-f f_*  f'_* \big\}
\cr&
\d(v+v_*-v'-v'_*)\,
\d\left({1\over 2}(v^2+v_*^2-v'\/^2-v_*'\/^2)\right).
}
\Eq(6.29)
$$

\bigskip

There remains to sum up the contributions of the terms
$\Cal W_3^\pi$ and
$\Cal V_3^\pi$.  It gives, after some computations using the exchange
of variables $v'\leftrightarrow v'_*$, the missing cross term
$$
\eqalign{&
(\Cal W^\pi_3+\Cal V_3^\pi)(x,v)=
2 \pi \int dv_*\int dv'\int dv_*' \;
[(f+f_*)f'f_*'-(f'+f_*')ff_*)]
\cr&
\qquad
\qquad
\d(v+v_*-v'-v_*') \,
\d\left({1\over2}(v^2+v_*^2-v'\/^2-v_*'\/^2)\right) \,
\widehat\phi(v'-v) \, \widehat\phi(v'-v_*).
}
$$
This completes the 
proof of the theorem.

\bigskip\bigskip


\heading 7. Concluding Remarks\endheading

\numsec=7
\reset
\numfor= 1
\bigskip


It is well known that other possible scaling lead to kinetic equations as 
well. The most
important is the low-density limit (or Boltzmann-Grad limit): it is the 
regime in which
classical rarefied gases are described by the usual Boltzmann equation.

In our
grand-canonical formalism it can be introduced in the following way.
We do not rescale
$\phi$ which is $O(1)$. On the other hand the rarefaction hypothesis is 
given by the
condition $\e^2\langle N\rangle =O(1)$. This means that (see \equ(2.24))
$$
f_0^\e(x,v)=O(\e).
\Eq(7.1)
$$
Under this scaling, the hierarchy becomes (see \equ(rp))
$$
\pt_tf^\e_j+\sum_{k=1}^jv_k\cdot \n_{x_k}f^\e_j=
{1\over \e}T_j^\e f_j^\e
+
{1\over\e^4}C^\e_{j+1}f^\e_{j+1}.
\Eq(7.2)
$$
Rescaling the correlation functions by defining:
$$
\bar f^\e_j=\e^{-j} f_j^\e,
\Eq(7.3)
$$
we arrive at the hierarchy
$$
\pt_t\bar f^\e_j+\sum_{k=1}^jv_k\cdot \n_{x_k}\bar f^\e_j=
{1\over \e}T_j^\e\bar f^\e_j
+
{1\over\e^3}C^\e_{j+1}\bar f^\e_{j+1},
\Eq(7.4)
$$
with a fixed initial datum of $O(1)$.
It is now clear that the terms $CC$ are vanishing in the limit $\e\to 0$ 
and the
statistical correlations are lost. On the other hand many terms of the 
type
$CT\dots T$ are finite in the limit.
It turns out that the sum of these terms lead to the Born series 
expansion of 
the cross section.
The underlying series actually converges
provided the potential $\phi$ is small.
This task is performed in the case of the Maxell-Boltzmann statistic 
in
Ref. [\rcite{BCEP2}] by the authors.
Here, a difficult point lies in the identification of the cross section as 
the Born series expansion of quantum scattering, a task which is achieved
using an original identity derived in [\rcite{Ca2}].

Another comment is in order.
The U-U equation has been partially derived whenever $f_0$
is the Wigner transform of a one-particle quasi-free state.
As shown in Appendix, a
sufficient condition for the explicit construction of such a state
is a small value of
the activity $z$. On the other hand the U-U equation for Bosons makes
sense also for more
general initial conditions describing states with large activity.
It seems very
interesting to understand whether the U-U dynamics of such states
make sense from a
physical point of view and whether it can describe dynamical
condensation phenomena.

\bigskip\bigskip


\heading {Appendix: Quasi-free states for Bosons}\endheading

\vskip .5cm
\resetall
\numfor= 1

\bigskip


Let $r$ be a one-particle state i.e. a self-adjoint positive operator
whose kernel is denoted by $r(x,y)$. We want to construct a state which is
compatible with the B-E statistics and with a given average particle
number.

Let $\s_{n}$ be a $n$-particle completely symmetric state given by
$$
\s_{n}(X_{n},Y_{n})=
\sum_{\pi\in \Cal P_{n}}
r(x_{1},y_{\pi(1)})\dots r(x_{n},y_{\pi(n)}).
\Eqa(A.1)
$$ 
The state
$$
\s^{z}=
{1\over \Xi(z)}\bigoplus_{n=0}^{\infty}{z^{n}\over n!}
\s_{n},
\Eqa(A.2)
$$
where 
$$
\Xi(z)=\sum_{n\ge 0}{z^n\over n!}\Tr \s_n,
\Eqa(A.3)
$$
is
a normalized state for Bosons and
$$
\langle N \rangle=\Tr [\s^z N]= z{d\over dz}\log \Xi.
\Eqa(A.4)
$$

We now compute the partition function $\Xi(z)$ (see [\rcite{G}]):
$$
\Xi(z)=\sum_{n\ge 0}{z^{n}\over n!}
\int dx_{1}\dots dx_{n}\;
\sum_{\pi\in \Cal P_{n}}
r(x_{1},x_{\pi(1)})\dots r(x_{n},x_{\pi(n)}).
\Eqa(A.5)
$$
Given $\pi$, let $\a_{1},\dots,\a_{n}$, be non negative integers,
$\a_{j}$ denoting the number of cycles of length $j$ in $\pi$.
Clearly
$$
\sum_{j=1}^{n}j\a_{j}=n.
\Eqa(A.6)
$$
Given the sequence $\a_{1},\dots,\a_{n}$,
$$
\int dx_{1}\dots dx_{n} \;
r(x_{1},x_{\pi(1)})\dots
r(x_{n},x_{\pi(n)})
=
\prod_{j=1}^{n}\Big(\Tr r^{j}\Big)^{\a_{j}}.
\Eqa(A.7)
$$
The number of permutations associated to a given sequence
$\a_{1},\dots,\a_{n}$ is
$$
n!\prod_{j=1}^{n}{1\over \a_{j}! j^{\a_{j}}}.
\Eqa(A.8)
$$
Hence
$$
\eqalign{
\Xi(z)=&
\sum_{n\ge 0}\sum\Sb \a_{1},\dots,\a_{n}\ge 0\\
\sum j\a_{j}=n\endSb \prod_{j=1}^{n}
{\Big(\Tr r^{j}\Big)^{\a_{j}}(z^{j})^{\a_{j}}\over \a_{j}! j^{\a_{j}}}
\cr&
=
\sum_{n\ge 0}\sum_{s\ge 0}{1\over s!}
\sum\Sb j_{1},\dots,j_{s}\ge 1\\
\sum j_{k}=n \endSb
\prod_{\ell=1}^{s}
{\Tr r^{j_{\ell}}z^{j_{\ell}}\over j_\ell}
=\exp\left[\sum_{j\ge 1}{\Tr r^{j}z^{j}\over j}\right].
}
\Eqa(A.9)
$$
In the second equality $s$ denotes the number of actual cycles in each 
permutation and $j_1,\dots,j_s$ are the lenghts of the cycles. The last 
sum is convergent for $z$ sufficiently small (away from Bose condensation 
region).
Then, by \equ(A.4) and \equ(A.9)
$$
\langle N\rangle=\sum_{j\ge 1}\Tr(r z)^j=z+o(z),
\Eqa(A.10)
$$
for $z$ small.

The RDM's according to \equ(2.5) are:
$$
\r_{j}(X_{j},Y_{j})=
{1\over \Xi(z)}\sum_{n\ge 0}{(n+j)!\over n!}
\int dZ_{n} \;
{z^{n+j}\over (n+j)!}
\s_{j+n}(X_{j},Z_{n};Y_{j},Z_{n}).
\Eqa(A.11)
$$
Therefore we have:
$$
\r_{j}(X_{j},Y_{j})=
{1\over \Xi(z)}
\sum_{n\ge j}{z^{n}\over (n-j)!}
\sum_{\pi\in \Cal P_{n}}
\int dz_{j+1}\dots \int dz_{n} \;
r(x_{1},\x_{\pi(1)})\dots r(z_{n},\x_{\pi(n)}),
\Eqa(A.12)
$$ 
where $\x=(Y_{j},Z_{n-j})$.
Hence
$$
\eqalign{
\r_{j}(X_{j},Y_{j})&=
{1\over \Xi(z)}\sum_{n\ge j}
{z^{n}\over (n-j)!}
\sum_{s=0}^{n-j}\left(\matrix n-j \\
s\endmatrix\right)
\sum_{\pi'\in\Cal P_{s}}
\int dZ_{s}\;
r(z_{1},z_{\pi'(1)})\dots r(z_{s},z_{\pi'(s)})
\cr&
(n-j-s)!
\sum_{\pi\in \Cal P_{j}}
\sum\Sb k_{1},\dots,\k_{j}\ge 1\\
\sum k_j=n-s\endSb\prod_{i=1}^{j} r^{k_{1}}(x_{i},y_{\pi(i)})
}
\Eqa(A.13)
$$
with $r^{k}(x,y)$ the kernel of $r^{k}$.

Since
$$
\left(\matrix n-j \\ s\endmatrix\right)
{(n-j-s)!\over (n-j)!}
=
{1\over s!},\quad z^{n}=z^{s}z^{n-s}=z^{s}
z^{\Sigma k_{\ell}},
\Eqa(A.14)
$$
we obtain
$$
\eqalign{
\r_{j}(X_{j},Y_{j})&=
{1\over \Xi(z)}
\sum_{s\ge 0}{z^{s}\over s!}
\Tr \s_{s}\sum_{\pi\in \Cal P_{j}}\prod_{\ell=1}^{j}
\sum_{k_{\ell}=1}^{\infty}z^{k_{\ell}}r^{k_{\ell}}(x_{\ell},y_{\pi(\ell)}).
}
\Eqa(A.15)
$$
Defining the one-particle operator
$$
r_{z}=\sum_{k\ge 1}z^{k}r^{k}={z r\over 1-z r},
\Eqa(A.16)
$$
for $z$ sufficiently small, we arrive at
$$
\r_{j}(X_{j},Y_{j})=
\sum_{\pi\in \Cal P_{j}}r_{z}(x_{1},y_{\pi(1)})
\dots r_{z}(x_{j}, y_{\pi(j)}),
\Eqa(A.17)
$$ 
that is the characterization of the quasi-free state in terms of
RDM.

\end

\bibitem[AM]{AM}
N.W. Ashcroft, N.D. Mermin,
{\bf Solid state physics}, Saunders, Philadelphia (1976).

\bibitem[Bo]{Bo}
A. Bohm, {\bf
Quantum Mechanics}, Texts and mo\-no\-graphs in Phy\-sics,
Sprin\-ger-Verlag (1979).

\bibitem[CC]{CC}
S. Chapman and T. G. Cowling, {\bf The Mathematical
Theory of
Non-uniform
Gases}, Cambridge Univ. Press, Cambridge, England (1970).

\bibitem[CIP]{CIP}
C. Cercignani, R. Illner, M. Pulvirenti,
{\bf The mathematical theory of dilute gases},
Applied Mathematical Sciences, Vol. 106,
Springer-Verlag, New York (1994).

\bibitem[Ch]{Ch}
S.L. Chuang, {\bf
Physics of optoelectronic
devices}, Wiley series in pure and applied optics,
New-York (1995).

\bibitem[CTDL]{CTDL}
C. Cohen-Tannoudji, B. Diu, F. Lalo\"e,
{\bf M\'ecanique Quantique},
I et II, Enseignement des Sciences, Vol. 16, Hermann (1973).

\bibitem[Co]{Co}
M. Combescot,
{\it
On the generalized golden rule for transition probabilities},
Phys. A: Math. Gene., Vol. 34, N. 31, pp. 6087-6104 (2001).

\noi
[D\"u]
R. D\"umcke,
{\it The low density limit for an $N$-level system
interacting with a free Bose or Fermi gas},
Comm. Math. Phys., Vol. 97, N. 3, pp. 331-359 (1985).

\bibitem[EY1]{EY1}
L. Erd\"os, H.T. Yau, {\it
Linear Boltzmann Equation as Scaling Limit of Quantum Lorentz Gas}
Advances in differential equations and mathematical physics (Atlanta,
GA, 1997), pp. 137-155, Contemp. Math., 217,
Amer. Math. Soc., Providence, RI (1998).

\bibitem[EY2]{EY2}
L. Erd\"os, H.T. Yau, {\it
Linear Boltzmann equation as the weak coupling limit of
a random Schr\"odinger equation},
Comm. Pure Appl. Math., Vol. 53, N. 6, pp. 667-735
(2000).

\bibitem[Fi]{Fi}
M.V. Fischetti, {\it Theory of electron transport
in small semiconductor devices using the Pauli master equations},
J. Appl. Phys., Vol. 83, N. 1, pp. 270-291 (1998).

\bibitem[HLW]{HLW}
T.G. Ho, L.J. Landau, A.J. Wilkins,
{\it On the weak coupling limit for a Fermi gas in a random potential},
Rev. Math. Phys., Vol. 5, N. 2, pp. 209-298 (1993).

\noi
[Hu]
K. Huang, {\bf
Statistical mechanics}, Wiley and Sons (1963).

\noi
[KL1]
W. Kohn, J.M. Luttinger, Phys. Rev., Vol. 108, pp. 590 (1957).

\noi
[KL2] W. Kohn, J.M. Luttinger,
Phys. Rev., Vol. 109, pp. 1892 (1958).

\noi
[Ku] R. Kubo, J. Phys. Soc. Jap., Vol. 12 (1958).

\bibitem[L]{L}
O. Lanford III,
{\bf Time evolution of large classical systems},
Lecture Notes in Physics, Vol. 38, pp. 1-111, E.J. Moser ed.,
Springer-Verlag (1975).

\noi
[La] L.J. Landau, {\it Observation of Quantum Particles on a
Large Space-Time Scale}, J. Stat. Phys., Vol. 77, N. 1-2, pp. 259-309 (1994).

\bibitem[MRS]{MRS}
P.A. Markowich, C. Ringhofer, C. Schmeiser,
{\bf Semiconductor equations},
Sprin\-ger-Verlag, Vienna  (1990).

\bibitem[Ni1]{Ni1}
[Ni1] F. Nier, {\it Asymptotic Analysis of a scaled Wigner equation and
Quantum Scattering}, Transp. Theor. Stat. Phys.,
Vol. 24, N. 4 et 5, pp. 591-629 (1995).

\bibitem[Ni2]{Ni2}
F. Nier, {\it A semi-classical picture of quantum scattering},
Ann. Sci. Ec. Norm. Sup., 4. S\'er., t. 29, p. 149-183 (1996).

\bibitem[RV]{RV}
E. Rosencher, B. Vinter, {\bf Optoelectronique},
Dunod (2002).

\bibitem[Sp1]{Sp1}
H. Spohn, {\it Derivation of the transport equation
for electrons moving through random impurities}, J. Stat. Phys.,
Vol. 17, N. 6, pp. 385-412 (1977).

\bibitem[Sp2]{Sp2}
H. Spohn, {\bf Large scale dynamics of interacting particles},
Springer (1991).

\bibitem[Sp3]{Sp3}
H. Spohn,
{\it Kinetic equations from Hamiltonian dynamics: Markovian limits},
Rev. Modern Phys., Vol.  52, N. 3, pp. 569-615 (1980).

\bibitem[Sp4]{Sp4}
H. Spohn,
{\it Quantum Kinetic Equations},  in ``On three levels: Micro-Meso and
Macro Approaches in Physics'', M. Fannes, C. Maes, A. Verbeure eds, Nato
ASI Series B: Physics Vol. 324, 1-10, (1994).

\noi
[VH1] L. Van Hove, Physica, Vol. 21 p. 517 (1955).

\noi
[VH2] L. Van Hove, Physica, Vol. 23 p. 441 (1957).

\noi
[VH3] L. Van Hove, in {\bf Fundamental Problems
in Statistical Mechanics}, E.G.D. Cohen ed., p. 157 (1962).

\endRefs
\end